\documentclass[12pt,a4paper,reqno]{amsart}

\usepackage[english]{babel}
\usepackage{amssymb,latexsym,amsfonts,amsthm,upref,amsmath}
\usepackage[foot]{amsaddr}
\usepackage[margin=1in]{geometry} 
\usepackage[dvipsnames,x11names]{xcolor}
 \definecolor{myblue}{HTML}{003399}
\usepackage{hyperref}
\hypersetup{colorlinks,citecolor=myblue,filecolor=black,linkcolor=myblue,urlcolor=myblue}
\usepackage{enumerate}  
\usepackage{tikz}
\usepackage{float}
\usepackage{cleveref}
\usepackage{mathtools}
\usepackage{cite}
\usepackage{filecontents}
\makeatletter
\newcommand{\leqnomode}{\tagsleft@true}
\newcommand{\reqnomode}{\tagsleft@false}

\makeatother
\newtheorem*{thm*}{Theorem}
\newtheorem*{lem*}{Lemma}
\newtheoremstyle{prim}{}{}{\normalfont}{}{\bfseries}{.}{ }{}
\newtheoremstyle{stil}{}{}{\slshape}{}{\bfseries}{.}{ }{}
\theoremstyle{stil}
\newtheorem{thm}{Theorem}[section]
\newtheoremstyle{defi}{}{}{}{}{\bfseries}{.}{ }{}
\theoremstyle{defi}
\newtheorem{defn}[thm]{Definition}
\theoremstyle{defi}
\newtheorem{rem}[thm]{Remark}
\theoremstyle{stil}
\newtheorem*{mthm*}{Main Theorem}
\newtheorem*{kor*}{Corollary}
\newtheorem{pro}[thm]{Proposition}
\theoremstyle{stil}
\newtheorem{lem}[thm]{Lemma}
\theoremstyle{stil}
\newtheorem{kor}[thm]{Corollary}
\theoremstyle{prim}

\newenvironment{prf}{\noindent \textit{Proof.}}{\null\hfill$\qed$\hskip
2mm\vskip 2mm}

\newcommand{\ur}{{\rm U}(R)}

\newcommand{\vr}{\mathcal{V}^c(R)}
\newcommand{\eps}{\varepsilon}
\newcommand{\iii}{\mathrm{i}}

\newcommand{\ahp}{\mathcal{A}_{h,p}(\widehat{\mathfrak{gl}}_2)}
\newcommand{\Uc}{ \mathcal{U} }

\newcommand{\Y}{{\rm Y}_h(\mathfrak{gl}_N)}

\newcommand{\Vc}{\mathcal{V}^{c} (R)}
 
\newcommand{\vac}{ \mathrm{\boldsymbol{1}}}

\newcommand{\gl}{\mathfrak{gl}}

\newcommand{\CC}{\mathbb{C}}
\newcommand{\RR}{\mathbb{R}}
\newcommand{\ZZ}{\mathbb{Z}}

\newcommand{\Lc}{\mathcal{L}}

\newcommand{\Sc}{\mathcal{S}}

\newcommand{\Ec}{\mathcal{E}}

\newcommand{\wht}{\widehat}
\newcommand{\wvr}{\overline}

\newcommand{\ot}{\otimes}
\newcommand{\ts}{\hspace{1pt}}

\newcommand{\tr}{ {\rm tr}}

\newcommand{\ndo}{\mathop{\mathrm{End}}}
\newcommand{\om}{\mathop{\mathrm{Hom}}}
\newcommand{\im}{\mathop{\mathrm{Im}}}

\newcommand{\diag}{\mathop{\mathrm{diag}}}

\newcommand{\cdotlr}{\mathop{\hspace{-2pt}\underset{\text{LR}}{\cdot}\hspace{-2pt}}}

\newcommand{\fand}{\quad\text{and}\quad}
\newcommand{\Fand}{\qquad\text{and}\qquad}

\newcommand{\non}{\nonumber}
\newcommand{\beq}{\begin{equation}}
\newcommand{\eeq}{\end{equation}}
\newcommand{\ben}{\begin{equation*}}
\newcommand{\een}{\end{equation*}}

\newcommand{\R}{\wvr{R}}
\newcommand{\Rv}{\check{R}}

\makeatletter
\def\smalloverbrace#1{\mathop{\vbox{\m@th\ialign{##\crcr\noalign{\kern3\p@}%
  \tiny\downbracefill\crcr\noalign{\kern3\p@\nointerlineskip}%
  $\hfil\displaystyle{#1}\hfil$\crcr}}}\limits}
\makeatother

\makeatletter
\def\smallunderbrace#1{\mathop{\vtop{\m@th\ialign{##\crcr
   $\hfil\displaystyle{#1}\hfil$\crcr
   \noalign{\kern3\p@\nointerlineskip}% 
   \tiny\upbracefill\crcr\noalign{\kern3\p@}}}}\limits}
\makeatother

\makeatletter
\def\author@andify{%
  \nxandlist {\unskip ,\penalty-1 \space\ignorespaces}%
    {\unskip {} \@@and~}%
    {\unskip \penalty-2 \space \@@and~}%
}
\makeatother

\begin{document}

\title{Associating modules for the $h$-Yangian and quantum elliptic algebra in type $A$ with $h$-adic quantum vertex algebras}

\author{Lucia Bagnoli}
\address[L. Bagnoli]{Dipartimento di Matematica, Sapienza Universit\`{a} di Roma, P.le Aldo Moro 5, 00185 Rome, Italy}
\email{lucia.bagnoli@uniroma1.it}

\author{Naihuan Jing}
\address[N. Jing]{Department of Mathematics, North Carolina State University,  Raleigh, NC 27695, USA}
\email{jing@ncsu.edu}
 
\author{Slaven Ko\v{z}i\'{c}}
\address[S. Ko\v{z}i\'{c}]{Department of Mathematics, Faculty of Science, University of Zagreb,  Bijeni\v{c}ka cesta 30, 10000 Zagreb, Croatia}
\email{slaven.kozic@math.hr}

\begin{abstract}
We consider the Etingof--Kazhdan quantum vertex algebra $\mathcal{V}^c(R)$ associated with the trigonometric and elliptic $R$-matrix of type $A.$ We establish a connection between (restricted) modules for the $h$-Yangian $\textrm{Y}_h(\mathfrak{gl}_N)$ and the elliptic quantum algebra $\mathcal{A}_{h,p}(\widehat{\mathfrak{gl}}_2)$  of level zero, and deformed (twisted) $\phi$-coordinated $\mathcal{V}^c(R)$-modules. As its application, in the trigonometric case,  we construct new families of central elements of  $\mathcal{V}^c(R)$ at the critical level $c=-N,$ which we then use to derive commutative families in the $h$-Yangian $\textrm{Y}_h(\mathfrak{gl}_N).$
\end{abstract}

\maketitle
\allowdisplaybreaks

%%%%%%%%%%%%%%%%%%%%%%%%%%%%%%%%%%%
%%%%%%%%%%%%%%%%%%%%%%%%%%%%%%%%%%%
\section{Introduction}\label{intro}
\setcounter{equation}{0}
\numberwithin{equation}{section}
%%%%%%%%%%%%%%%%%%%%%%%%%%%%%%%%%%%
%%%%%%%%%%%%%%%%%%%%%%%%%%%%%%%%%%%
The definition and first examples of quantum vertex algebras, the {\em quantum affine vertex algebras}, were introduced by Etingof and Kazhdan \cite{EK}. 
Later on, Li \cite{Liphi} developed the theory of {\em $\phi$-coordinated modules},
which was instrumental in associating quantum vertex algebra theory with representations of quantum affine algebras; see  \cite{JKLT,Kong,K} and references therein.
The goal of this paper is to establish a connection between the 
aforementioned Etingof--Kazhdan  quantum vertex algebras $\vr,$ associated with the trigonometric and elliptic $R$-matrix $R$   of type $A,$
and representation theories of
  the $h$-Yangian  $\Y$ and the level zero quantum elliptic algebra $\ahp .$   

Both classes of these quantum algebras are characterized by {\em multiplicative} defining relations, which can be written in the {\em FRT-form} \cite{FRT}   using the  $R$-matrix formalism. On the other hand,    $\vr$ is essentially   {\em additive} in nature  and its $\mathcal{S}$-locality property  
  takes the form of the {\em quantum current commutation relation} of Reshetikhin and Semenov-Tian-Shansky \cite{RS}. To reconcile these fundamental differences, we employ the notion of {\em deformed $\phi$-coordinated module} \cite{BK}. Roughly speaking,  we use the $\phi$-coordinated module  theory to connect  the additive and multiplicative setting, and the   notion of deformed module,  which was motivated in part by the work of Anguelova and Bergvelt \cite{AB}, to establish a relation between     the FRT-relations  and the quantum current commutation relations.
Furthermore, in the elliptic case, due to the   form of the underlying  $R$-matrix of the eight-vertex model, we employ the notion of
{\em twisted} deformed $\phi$-coordinated module, which is also motivated by the    twisted  $\phi$-coordinated modules of Li, Tan and Wang \cite{LiTW}, and we slightly extend the original Etingof--Kazhdan construction of $\vr .$

The main results of this paper are Theorems \ref{thm_53-trigg}, \ref{thm_53} and
\ref{thm_66}, which establish a connection between (twisted) deformed $\phi$-coordinated $\vr$-modules and (restricted) modules for the algebras $\Y$ and   $\ahp$.
Moreover, in the trigonometric case, we demonstrate an application of these results to constructing generators of the quantum {\em  Feigin--Frenkel center} \cite{FF}, i.e. the center of $\vr$ at the critical level $c=-N,$ and commutative families in the $h$-Yangian $\Y.$ Finally, the fusion procedure of  Frappat, Issing and Ragoucy \cite{FIR} indicates a possibility of similar applications, as well as  generalizations of our main results, to the case of the elliptic algebra $\mathcal{A}_{q,p}(\widehat{\mathfrak{gl}}_N)$.

%%%%%%%%%%%%%%%%%%%%%%%%%%%%%%%%%%%
%%%%%%%%%%%%%%%%%%%%%%%%%%%%%%%%%%%
\section{Preliminaries}\label{sec01}
%%%%%%%%%%%%%%%%%%%%%%%%%%%%%%%%%%%
%%%%%%%%%%%%%%%%%%%%%%%%%%%%%%%%%%%

In this section, we recall the trigonometric and elliptic $R$-matrices of type $A.$
We adapt their definitions so that they are given over the commutative ring $\CC[[h]]$ of formal power series in parameter $h$, instead of over the field  $\CC(q^{1/2})$ of rational functions in $q^{1/2}.$   In both cases, the parameters $h$ and $q$ are related via $q^{1/2}=e^{h/4}.$ Throughout the paper, in the trigonometric case, we denote by $N$  a fixed integer greater than or equal to $2,$ while, in the elliptic case, we set $N=2.$

%%%%%%%%%%%%%%%%%%%%%%%%%%%%%%%%%%%
\subsection{Trigonometric \texorpdfstring{$R$}{R}-matrix}\label{sec0101}
%%%%%%%%%%%%%%%%%%%%%%%%%%%%%%%%%%%

In this subsection, we recall the trigonometric $R$-matrix of   type $A$, which appears in the FRT-realization of the quantum affine algebra $U_q(\widehat{\mathfrak{gl}}_N)$; for more information see the references \cite{FRT,J,PS,RS}. 

By \cite{FR}, there exists a unique formal power series 
\beq\label{efodze00}
f_q(z)=1+\sum_{r=1}^\infty f_{q,r} \frac{z^r}{(1-z)^{r}}\in\CC(q)[[z]],
\eeq  
  such that all $f_{q,r}(q-1)^{-r}$ are regular at $q=1$, satisfying
$$
f_q(zq^{2N})=f_q(z) (1-zq^2)(1-zq^{2N-2}) (1-z)^{-1}(1-zq^{2N})^{-1}.
$$
By setting $q=e^{h/2},$ we obtain  
$$
f(z)\coloneqq f_{e^{h/2}}(z)=1+\sum_{r=1}^\infty f_{r} \frac{z^r}{(1-z)^{r}}\in\CC[[h,z]] ,\quad\text{where}\quad
f_r=f_{q,r}\big|_{q=e^{h/2}}.
$$
 For more information on the series $f(z)$ see also \cite{KM}.

 Let $V=\ndo\CC^N\ot\ndo\CC^N.$ 
Consider  the   trigonometric $R$-matrix   $\R (z )\in V  [[z,h]],$    
\begin{align}
\R (z ) =&\,\sum_{i=1}^N e_{ii}\ot e_{ii} 
+ e^{-h/2}\frac{1-z}{1-ze^{-h}}\sum_{\substack{i,j=1\\i\neq j}}^N e_{ii}\ot e_{jj} \non\\
&+ \frac{\left(1-e^{-h }\right)z }{1-ze^{-h}} \sum_{\substack{i,j=1\\i> j}}^N e_{ij}\ot e_{ji}
+\frac{ 1-e^{-h } }{1-ze^{-h}} \sum_{\substack{i,j=1\\i< j}}^N e_{ij}\ot e_{ji},\label{rbar_trigg}
\end{align}
where $e_{ij}\in\ndo\CC^N$ are the matrix units. 
Finally,  define the normalized $R$-matrix by
\beq\label{erofzed}
R(z)=  f(z)   \ts \R(z) .
\eeq
Due to the form of the term  $f(z)$, the $R$-matrix $R(z)$ can be also  regarded as an element  of $V(z)[[h]].$ It satisfies
the (multiplicative) {\em quantum Yang--Baxter equation},  
\beq\label{trig_ybe_12}
R_{12}(z_1)\ts R_{13}(z_1 z_2)\ts R_{23}(z_2)
 =   R_{23}(z_2)\ts R_{13}(z_1 z_2)\ts R_{12}(z_1) 
\eeq
and the {\em crossing symmetry properties,}
\beq\label{csym}
R(ze^{ Nh})^{t_1} \ts D_1\ts  ( R(z)^{-1})^{t_1}=D_1\fand (R(z)^{-1})^{t_2} \ts D_2 \ts R(z e^{ Nh})^{t_2} = D_2.
\eeq
Here $t_k$ stands for the matrix transposition $e_{ij}\mapsto e_{ji}$ applied on the $k$-th tensor factor, 
$D_1=D\ot 1, $   $D_2=1\ot D $
and $D$ is the diagonal $N\times N$ matrix
\beq\label{dijagonalna}
D=\diag\left(e^{ (N-1)h/2 },e^{ (N-3)h/2},\ldots ,e^{- (N-1)h/2} \right) .
\eeq

Consider the unique embedding 
$\CC_*(u)\hookrightarrow\CC((u))  ,$ where $\CC_*(u)$ denotes   the localization of the ring of formal Taylor series $\CC[[u]]$ at $\CC[u]\setminus\left\{0\right\}.$ 
It naturally extends to the embedding
\beq\label{iotau_emb}
\iota_u \colon \CC_*(u)[[h]]\hookrightarrow\CC((u)) [[h]].
\eeq
Let us regard the expression in \eqref{erofzed} as a formal Taylor   series in $h.$
By applying the substitution $z=e^u$ and then the map $\iota_u ,$ one obtains an element $R^\prime(u)$ of
$V ((u))[[h]]$;  see \cite{KM} for more details. 
Furthermore, by \cite[Prop. 1.2]{EK4} and \cite[Prop. 2.1]{KM}, there exists a unique $\psi\in 1+h\CC[[h]]$ such that the $R$-matrix 
\beq\label{erofu}
R(e^u)\coloneqq \psi \ts R^\prime(u)\in \ndo\CC^N \ot\ndo\CC^N ((u))[[h]]
\eeq
has the {\em unitarity property},
$
R_{12}(e^u) R_{21}(e^{-u}) =1,
$
and 
the {\em crossing symmetry properties},
\beq\label{csym_reu}
R(e^{u+Nh})^{t_1} \ts D_1\ts  ( R(e^u)^{-1})^{t_1}=D_1\fand (R(e^u)^{-1})^{t_2} \ts D_2 \ts R( e^{u+Nh})^{t_2} = D_2.
\eeq
In addition, the $R$-matrix \eqref{erofu}   
satisfies the  (additive) {\em quantum Yang--Baxter equation},
$$
R_{12}(e^u)\ts  R_{13}(e^{u+v}) \ts R_{23}(e^v)=R_{23}(e^v)\ts R_{13}(e^{u+v}) \ts R_{12}(e^u).
$$

%%%%%%%%%%%%%%%%%%%%%%%%%%%%%%%%%%%
\subsection{Elliptic \texorpdfstring{$R$}{R}-matrix}\label{sec0102}
%%%%%%%%%%%%%%%%%%%%%%%%%%%%%%%%%%%

 In this section, we follow  Foda et al. \cite[Sect. 2]{FIJKMY} to 
 recall the elliptic $R$-matrix of the 
 eight-vertex
  model; for more information see  \cite{Bax,Bel,Chu, RT,T}.

In order to introduce the $R$-matrix, we   need the {\em infinite $q$-Pochhammer symbols},
$$
(z;p )_\infty
=\prod_{k\geqslant 0} (1-zp^{k} ) \in \CC[[z,p]].
$$
Consider the formal power series, which is  obtained from \eqref{efodze00} for $N=2,$ 
$$
f(z)=\frac{(z;q^{4})_\infty\ts (zq^{4};q^{4})_\infty}{(zq^{2};q^{4})_\infty^2}\in\CC(q)[[z]].
$$
Let $p$ be another formal parameter,   the {\em elliptic nome}.
Introduce the formal  series
$$
\alpha(z) = \frac{(p^{1/2}qz;p)_\infty}{(p^{1/2}q^{-1}z;p)_\infty}\prod_{k\geqslant 1} f(p^k z^2)\Fand
\beta(z)= \frac{(-pqz;p)_\infty}{(-pq^{-1}z;p)_\infty}
\prod_{k\geqslant 1} f(p^k z^2).
$$
Clearly, $\alpha(z)$ and $\beta(z)$ belong to $ \CC(q )[z ][[p^{1/2}]].$
Finally, the {\em elliptic $R$-matrix}  is given, with respect to the standard basis $e_1\ot e_1,e_1\ot e_2,e_2\ot e_1,e_2\ot e_2$ of the space $\CC^2\ot\CC^2$, where $e_1=(1,0)$ and $e_2=(0,1),$ by
\beq\label{R}
R(z)=\begin{pmatrix}
a(z)& 0 & 0 & d(z)\\
 0 & b(z) & c(z) & 0\\
 0 & c(z) & b(z) & 0\\
d(z)& 0 & 0 & a(z)
\end{pmatrix}.
\eeq
Its matrix entries  are  uniquely determined by the identities
$$
a(z)+d(z)=q^{-1/2}f(z^2)^{-1}\frac{\alpha(z^{-1})}{\alpha(z)}\fand
b(z)+c(z)=q^{1/2}f(z^2)^{-1} \frac{1+q^{-1}z}{1+qz}\frac{\beta(z^{-1})}{\beta(z)} 
$$
and a  requirement that  $a(z)$ and $b(z)$ (resp. $c(z)$ and $d(z)$) have only even (resp. odd) powers of $z.$  All entries belong to $\CC(q^{1/2})((z))[[p^{1/2}]]$ and, furthermore, they can be expressed in terms of the Jacobi theta function; see  \cite[Sect. 2.1]{FIJKMY2} for the explicit formulas.

 We   now   adjust the above setting, so that   the   $R$-matrix \eqref{R} is defined over the  ring $\CC[[h]]$; we   omit some technical details as they are  available in \cite[Sect. 2.2]{BBK}. First of all, let $q^{1/2}=e^{h/4}\in\CC[[h]].$ Next, we fix $a\in\RR_{>0}$ and $b\in\ZZ_{>0}$ and then set 
\beq\label{pandh}
p=a h^{2b}\in\CC[[h]].
\eeq
The $R$-matrix \eqref{R} is now defined over $\CC[[h]] $  and   we have
\beq\label{Rzh}
R(z)\in\ndo\CC^2 \ot  \ndo\CC^2  ((z ))[[h]].
\eeq
It satisfies
the (multiplicative) {\em quantum Yang--Baxter equation},  
\beq\label{ell_ybe_12}
R_{12}(z_1)\ts R_{13}(z_1 z_2)\ts R_{23}(z_2)
 =   R_{23}(z_2)\ts R_{13}(z_1 z_2)\ts R_{12}(z_1),
\eeq
and it possesses the   property 
\beq\label{uniz}
R_{12}(z)\ts R_{21}(1/z)=\Uc(z), \qquad\text{where}\qquad \Uc(z)=  
 e^{-h/2} f(z^2)^{-1} f(z^{-2})^{-1}.
\eeq

We shall also need an additive counterpart  of  $R(z)$. It is   obtained by setting $z=(-1)^\eps e^u\in\CC[[u]],$  $\eps=0,1,$ in \eqref{Rzh}. The resulting $R$-matrix $R(e^{u+\eps\pi\iii})\coloneqq R((-1)^\eps e^u),$ where $\iii\in\CC$ denotes the imaginary unit,  satisfies
\beq\label{Ruh}
R(e^{u+\eps\pi\iii})\in\ndo\CC^2 \ot  \ndo\CC^2  ((u))[[h]].
\eeq
Moreover, it satisfies
the (additive) {\em quantum Yang--Baxter equation},  
\begin{align}
&R_{12}(e^{u_1+\eps_1\pi \iii})\ts R_{13}(e^{u_1+u_2+(\eps_1+\eps_2)\pi \iii})\ts R_{23}(e^{u_2+\eps_2\pi \iii})\non\\
 =&\,  R_{23}(e^{u_2+\eps_2\pi \iii})\ts R_{13}(e^{u_1+u_2+(\eps_1+\eps_2)\pi \iii})\ts R_{12}(e^{u_1+\eps_1\pi \iii}),\label{ell_addi}
\end{align}
where $\eps_1,\eps_2\in\left\{0,1\right\},$
and it possesses the   property 
$$
R_{12}(e^{u+\eps\pi\iii})\ts R_{21}(e^{-(u+\eps\pi\iii)})= \Uc(e^u), \qquad\text{where}\qquad \Uc(e^u)=  
 e^{-h/2} f(e^{2u})^{-1} f(e^{-2u})^{-1}.
$$ 

 \begin{rem}\label{trig_ell_rem}
Consider the $R$-matrix    \eqref{Rzh}. Due to \eqref{pandh}, the elliptic limit $p\to 0$ of $R(\pm z)$ is found by setting $a\to 0$. It produces the
trigonometric $R$-matrix
\begin{align*}
R^{trig}_{\pm}(  z)\coloneqq
\frac{e^{-h/4}}{f(z^2)}\left(
\sum_{i=1,2} e_{ii}\ot e_{ii}
+r (z)\sum_{\substack{i,j=1,2\\i\neq j}}e_{ii}\ot e_{jj}
\pm s  (z)\sum_{\substack{i,j=1,2\\i\neq j}}e_{ij}\ot e_{ji}
\right),
\end{align*} 
where
$$
r (z)=\frac{e^{h/2 }(1-z^2)}{(1-e^{h} z^2)}\fand s  (z)=   \frac{(1-e^{h})z}{(1-e^{h} z^2)},
$$
which, aside from the rescaled parameters, is of the same form as the  trigonometric $R$-matrix 
\eqref{erofzed} for $N=2$; see, e.g., \cite{AAFRR,FIR} for more information on   the connections between   trigonometric and elliptic settings. Next, by setting $z=e^u$ in $R^{trig}_{\pm}(  z)$ and then extracting the top degree components with respect to the degree operator given by
$\deg u^k h^l =-k-l,$
one obtains the rational $R$-matrix
$$
R^{rat}_\pm (u)\coloneqq\frac{1}{1+\frac{h}{ 2u}}\left(I+\frac{h}{2 u}P_\pm\right),
$$
where $P_+$ is the permutation operator,
\beq\label{p_plusminus}
P_+=\sum_{i,j=1,2} e_{ij}\ot e_{ji},\Fand P_-=\sum_{i,j=1,2} (-1)^{i+j} e_{ij}\ot e_{ji}.
\eeq
 Clearly, $R^{rat}_+ (u)$ is the  (normalized) Yang $R$-matrix. On the other hand, the $R$-matrix $R^{rat}_- (u)$ can be obtained by applying the Baxterization procedure \cite[Prop. 12]{GS} on the skew-invertible involutive symmetry $P_- .$
\end{rem}

%%%%%%%%%%%%%%%%%%%%%%%%%%%%%%%%%%%
\subsection{Notation}\label{sec0103}
%%%%%%%%%%%%%%%%%%%%%%%%%%%%%%%%%%%

From now on, we denote by $R^{trig}(z)=R(z)$ (resp. $R^{ell}(z)=R(z)$) the   trigonometric (resp. elliptic) $R$-matrix     \eqref{erofzed} (resp. \eqref{Rzh}), and by $R^{trig}(e^u)=R (e^u)$ (resp. $R^{ell}(e^{u+\varepsilon\pi\iii})=R (e^{u+\varepsilon\pi\iii})$) the   trigonometric (resp. elliptic) $R$-matrix     \eqref{erofu} (resp. \eqref{Ruh}).
To simplify the notation, we  omit the superscripts ``$trig$'' and ``$ell$'' whenever it is clear from the context whether the trigonometric or elliptic case is considered.
In the rest of this subsection, we introduce the notation which applies to both  settings.

For any positive integers $n$ and $m,$  let   
\beq\label{uovi}
u=(u_1,\ldots ,u_n)\fand v=(v_1,\ldots ,v_m)
\eeq
be the families of variables
 and
\beq\label{epsovi}
\eps=(\eps_1,\ldots ,\eps_n)\in\left\{0,1\right\}^{\times n}\fand 
\nu=(\nu_1,\ldots ,\nu_m)\in\left\{0,1\right\}^{\times m} 
\eeq
the binary tuples.
Define
\beq\label{uepsovi}
u^\eps =(u_1+\eps_1   \pi\iii,\ldots ,u_n+\eps_n   \pi\iii)\fand
v^\nu =(v_1+\nu_1   \pi\iii,\ldots ,v_m+\nu_m   \pi\iii).
\eeq
In the trigonometric setting, we  only consider the case $\eps=(0,\ldots ,0)$ and $\nu=(0,\ldots ,0),$ so the families \eqref{uovi} and \eqref{uepsovi} coincide, i.e., we have $u=u^\eps$ and $v=v^\nu . $
Let $z$ be a single variable.
We associate with the $(n+m)$-tuple $(u^\eps, v^\nu)$ and $d\in\CC$ the $R$-matrix product with coefficients in the tensor algebra $(\ndo\CC^N)^{\ot n} \ot(\ndo\CC^N)^{\ot m}$ (with $N=2$ in the elliptic case),
\beq\label{rnm-add}
R_{nm}^{12}(e^{z+u^\eps-v^\nu +dh})=
\prod_{p=1,\dots,n }^{\longrightarrow} \prod_{q=n+1,\dots,n+m }^{\longleftarrow}
R_{pq}(e^{z+ u_p-v_{q-n} +(\eps_p-\nu_{q-n})   \pi\iii +dh}),
\eeq
where the arrows indicate the order of the factors and the superscript $1$ (resp. $2$) corresponds to the tensor factors $1,\ldots ,n$ (resp. $n+1,\ldots ,n+m$).
Also, we write
\beq\label{rnm-addd}
R_{nm}^{12}(e^{ u^\eps-v^\nu +dh})=
\prod_{p=1,\dots,n }^{\longrightarrow} \prod_{q=n+1,\dots,n+m }^{\longleftarrow}
R_{pq}(e^{  u_p-v_{q-n} +(\eps_p-\nu_{q-n})   \pi\iii +dh}).
\eeq
Finally, we   extend this notation    to the case of multiplicative $R$-matrix
 by
\beq\label{rnm-adddd}
 R_{nm}^{12}(ze^{ u^\eps-v^\nu +dh})=
\prod_{p=1,\dots,n }^{\longrightarrow} \prod_{q=n+1,\dots,n+m }^{\longleftarrow}
R_{pq}(ze^{  u_p-v_{q-n} +(\eps_p-\nu_{q-n})   \pi\iii +dh}).
\eeq
In addition, for $x=(x_1,\ldots ,x_n)$ and $y=(y_1,\ldots ,y_m)$ we write
\beq\label{rnm-addddd}
 R_{nm}^{12}(x/y)=
\prod_{p=1,\dots,n }^{\longrightarrow} \prod_{q=n+1,\dots,n+m }^{\longleftarrow}
R_{pq}(x_p/y_{q-n}).
\eeq
For example, set $n=3$ and $m=2.$ Denote by $R_{pq}$ the factor corresponding to the
pair $(p,q)\in\left\{1,2,3\right\}\times\left\{4,5 \right\}$ on the right-hand side of
\eqref{rnm-add}--\eqref{rnm-addddd}. Then the   product defined by \eqref{rnm-add}--\eqref{rnm-addddd} takes the form
$R_{15}R_{14}R_{25}R_{24}R_{35}R_{34}.$

%%%%%%%%%%%%%%%%%%%%%%%%%%%%%%%%%%%
\section{Quantum algebras associated with   \texorpdfstring{$R$}{R}-matrices}\label{sec02}
%%%%%%%%%%%%%%%%%%%%%%%%%%%%%%%%%%%

In this section, we employ the  $R$-matrices from Section \ref{sec01} to  recall the $R$-matrix realizations of the 
{\em $h$-Yangian} $\Y$ and the
{\em elliptic quantum algebra}   $\ahp.$ 
In contrast with their original definitions, we consider  both algebras   over the commutative ring $\CC[[h]] ,$  so that they are in tune with the $h$-adic quantum vertex algebra theory    \cite{EK}.

%%%%%%%%%%%%%%%%%%%%%%%%%%%%%%%%%%%
\subsection{\texorpdfstring{$h$}{h}-Yangian \texorpdfstring{$\Y$}{Yh(glN)}}\label{sec0201}
%%%%%%%%%%%%%%%%%%%%%%%%%%%%%%%%%%%

To recall the definition of  the {\em $h$-Yangian} $\Y ,$   we make use of  the FRT-realization of the quantum affine algebra in type $A$, which goes back to Reshetikhin and  Semenov-Tian-Shansky \cite{RS}; see also the paper by I. Frenkel and Reshetikhin \cite{FR}.
The algebra $\Y$ is generated by the elements $\ell_{ij}^{(r)},$ 
where $i,j=1,\ldots ,N$ and $r=0,1,\ldots .$ They are organized into matrix of formal power series
\beq\label{elminusihY}
L^-(z)=\sum_{i,j=1}^N e_{ij}\ot \ell_{ij}(x),\quad\text{where}\quad
\ell_{ij}(z)=\delta_{ij}-h\sum_{r\geqslant 0} \ell_{ij}^{(r)}\ts z^{-r}.
\eeq
The defining relations\footnote{The original definition of the $h$-Yangian   contains additional relations $\ell_{ij}^{(0)}=\delta_{ij}$ for $1\leqslant i\leqslant j\leqslant N,$ which we omit as they are not needed in the setting of this paper. However, note that all diagonal elements $\ell_{ii}^{(0)}$ are   invertible if we assume that the $h$-Yangian is $h$-adically completed.} for $\Y$ are given in terms of the trigonometric $R$-matrix \eqref{erofzed}, 
\beq\label{rll_def_trig}
R(z_1/z_2)\ts L^-_1(z_1)\ts L^-_2(z_2)
=L^-_2(z_2)\ts L^-_1(z_1)\ts R(z_1/z_2).
\eeq
Throughout the paper, we often use the subscripts to indicate the     factors in the tensor product algebras such as $(\ndo\CC^N)^{\ot n}\ot \Y ,$  e.g., for any     $n\geqslant k> 0 ,$
we write
\beq\label{notconv}
L^-_k(z)=\sum_{i,j=1 }^N 1^{\ot (k-1)   }\ot e_{ij}\ot 1^{\ot (n-k)}\ot\ell_{ij}(z).
\eeq
In particular, the above notation is used in \eqref{rll_def_trig} with $n=2$ and $k=1,2.$

Replacing the $R$-matrix $R(z)$ in   \eqref{rll_def_trig} by $\wvr{R}(z) ,$   defined by \eqref{rbar_trigg}, one obtains equivalent relations,
as the series $L^-(z)$ possesses only nonnegative powers of  $z,$  so the normalization factors cancel.
Throughout the paper, the $h$-Yangian  is  assumed to be topologically free. 

Using the $R$-matrix notation \eqref{rnm-addddd}, the  relation \eqref{rll_def_trig} can be   generalized as follows. For any $m,n=1,2,\ldots$ and the  variables $z=(z_1,\ldots ,z_n)  $ and $w=(w_1,\ldots ,w_m)$ we have
\beq\label{hY_gen_rtt}
R_{nm}^{12}(z/w)\ts L_{[n]}^{-\ts 13}(z)\ts L_{[m]}^{-\ts 23}(w)=
L_{[m]}^{-\ts 23}(w)\ts L_{[n]}^{-\ts 13}(z)\ts R_{nm}^{12}(z/w),
\eeq
where
\begin{align}
& L_{[n]}^{-\ts 13}(z)= L^{-}_{1\ts n+m+1}(z_1)\ldots L^-_{n\ts n+m+1}(z_n),\label{superscripts1}\\
& L_{[m]}^{-\ts 23}(w) = L^{-}_{n+1\ts n+m+1}(w_1)\ldots L^{-}_{n+m\ts n+m+1}(w_m).\label{superscripts2}
\end{align}

%%%%%%%%%%%%%%%%%%%%%%%%%%%%%%%%%%%
\subsection{Elliptic quantum algebra   \texorpdfstring{$\ahp$}{Ahp(gl2)}}\label{sec0202}
%%%%%%%%%%%%%%%%%%%%%%%%%%%%%%%%%%%

In this subsection, we follow  Foda et al. 
\cite{FIJKMY,FIJKMY2} 
to
recall    the 
{\em elliptic quantum algebra} $\ahp.$
 As before, we assume that the parameters $h$ and $p$ are related by \eqref{pandh}.
For simplicity's sake,  we only give the definition at the level zero, which we need later on. The algebra $\ahp $ is generated by the elements 
$\bar{L}_{i j,n},$  where $i,j\in\left\{1,2\right\} $ and $n\in\mathbb{Z},$ 
which are organized into     power series,
\beq\label{Lijn}
L_{ij }(z)=\delta_{ij} - \sum_{n\in\ZZ} L_{ij,n}\ts z^{-n},\quad\text{where}\quad
L_{ij,n} = (-p^{1/2})^{\max\left\{n,0\right\}}\bar{L}_{ij,n}.
\eeq
In addition, we assume that
\beq\label{Lijnn}
\bar{L}_{ij,n}=0\quad\text{if}\quad (-1)^{i+j}\neq (-1)^n.
\eeq
The defining relations for   $\ahp$ are expressed in terms of the  matrix
$$
L(z)=\sum_{i,j=1,2}e_{ij}\ot L_{i j }(z) 
$$
and the elliptic $R$-matrix \eqref{Rzh}, using the tensor product  notation convention  \eqref{notconv}, as
\beq\label{rll_def}
R(z_1/z_2)\ts L_1(z_1)\ts L_2(z_2)
=L_2(z_2)\ts L_1(z_1)\ts R(z_1/z_2).
\eeq
Throughout the paper, the algebra $\ahp$ is assumed to be topologically free.

Suppose $U$ is a topologically free $\CC[[h]]$-module, i.e., equivalently, $U=U^0[[h]]$ for some complex vector space $U^0.$  We denote by $U((z_1,\ldots ,z_n))_h$ the $h$-adic completion of the $\CC[[h]]$-module of truncated formal Laurent series $U((z_1,\ldots ,z_n)) ,$ so that we have 
$U((z_1,\ldots ,z_n))_h=U^0((z_1,\ldots ,z_n))[[h]].$
In the rest of this section, we derive two simple properties of the generator matrix $L(z).$ First, we observe that, due to \eqref{Lijn}, for any $n=1,2,\ldots$ and the family of variables $z=(z_1,\ldots ,z_n) ,$ we have
\beq\label{elen}
L_1(z_1)\ldots L_n(z_n)\in\left(\ndo\CC^2\right)^{\ot n}\ot\om  \big(\ahp,\ahp((z_1,\ldots ,z_n))_h \big).
\eeq
We shall denote the expression in \eqref{elen} more briefly by 
$
L_{[n]}(z)=L_{[n]}(z_1,\ldots ,z_n) .$
Let $m$ be a positive integer and $w=(w_1,\ldots ,w_m)$ another family of variables.
As with the $h$-Yangian (recall \eqref{hY_gen_rtt}),
the defining relation   \eqref{rll_def} for $\ahp$ can be generalized as  
$$
R_{nm}^{12}(z/w)\ts L_{[n]}^{13}(z)\ts L_{[m]}^{23}(w)=
L_{[m]}^{23}(w)\ts L_{[n]}^{13}(z)\ts R_{nm}^{12}(z/w),
$$
where the meaning of the superscripts $1,2,3$ is the same as in \eqref{superscripts1} and \eqref{superscripts2}.

%%%%%%%%%%%%%%%%%%%%%%%%%%%%%%%%%%%
%%%%%%%%%%%%%%%%%%%%%%%%%%%%%%%%%%%
\section{Etingof--Kazhdan's quantum vertex algebra construction}\label{sec03}
%%%%%%%%%%%%%%%%%%%%%%%%%%%%%%%%%%%
%%%%%%%%%%%%%%%%%%%%%%%%%%%%%%%%%%%

In this section, we recall the Etingof--Kazhdan construction \cite{EK} of the quantum affine vertex algebra associated with  the trigonometric $R$-matrix of  type $A$. Next, in the elliptic case, we extend the original construction from  \cite{EK} to obtain an $h$-adic quantum vertex algebra which is  related with the elliptic quantum algebra  $\ahp,$ in the sense of  Theorems  \ref{thm_53} and \ref{thm_66} below. For a precise definition of the notion of ($h$-adic) quantum vertex algebra see \cite[Sect. 1.4.1]{EK} and \cite[Def. 2.20]{Li}. From now on, the tensor products of topologically free $\CC[[h]]$-modules are assumed to be $h$-adically completed, e.g., for $U=U^0[[h]],$ where $U^0$ is a complex vector  space, we denote $\left(U^0\ot_\CC U^0\right)[[h]]$ by $U\ot U.$

%%%%%%%%%%%%%%%%%%%%%%%%%%%%%%%%%%%
\subsection{Trigonometric case}\label{sec0301}
%%%%%%%%%%%%%%%%%%%%%%%%%%%%%%%%%%%

In this subsection, we consider the trigonometric $R$-matrix    \eqref{erofu}.
First, following \cite{EK3,EK4}, we introduce a certain   algebra $\textrm{U}(R)=\textrm{U}(R^{trig})$; see also \cite{FRT,RS}.
It is defined as an associative algebra over the ring $\CC[[h]]$ generated by the elements $t_{ij}^{(-r)}$,   $i,j=1,\ldots ,N$ and $r=1,2,\ldots ,$ subject to the  defining  relations
\beq\label{rtt}
R(e^{u-v})\ts T^+_{1} (u)\ts T^+_2  (v)=  T^+_2  (v)\ts T^+_{1} (u)\ts R(e^{u-v}),
\eeq
where the matrix $T^+ (u) $ is given by
\beq\label{tplustrig}
T^+ (u) =\sum_{i,j=1}^N e_{ij}\ot t^+_{ij}  (u)
\qquad\text{for}\qquad
 t^+_{ij} (u)=\delta_{ij}-h\sum_{r=1}^{\infty}t_{ij}^{(-r)}u^{r-1} .
\eeq
Clearly, $\textrm{U}(R)$ can be regarded as a trigonometric counterpart of the dual Yangian for the general linear Lie algebra $\gl_N$; see, e.g., \cite{I}. 
 
Denote by $\vac$  the unit in   $\textrm{U}(R).$    
Let $\Vc=\mathcal{V}^c(R^{trig}) $ be the $h$-adic completion
 of the $\CC[[h]]$-module of $\ur ,$  where the complex parameter $c $ in the superscript determines the action of the operator series $T^-(u),$ as  given by the next lemma  which goes back to  \cite[Lemma 2.1]{EK}. From now on, we  consider   $T^+ (u)$  as   operator series  over $\Vc,$
such that its action is given by the algebra multiplication.

\begin{lem}\label{lemma21}
For any $c\in\CC$ there exists a unique invertible operator series
$$T^- (u)\in\ndo\CC^N \ot \om( \Vc ,\Vc  ((u))_h )$$
such that   for all $n\geqslant 0$ we have
\beq\label{teminus}
 T^-_{0} (u)\ts T_{1}^+ (v_1)\ldots T_{n}^+ (v_n)\vac=
R^{01}_{1n}(e^{u-v+hc/2})^{-1}
\ts T_{1}^+ (v_1)\ldots T_{n }^+ (v_n) \ts  R^{01}_{1n}(e^{u-v-hc/2})\vac.
\eeq
\end{lem}

For any $n\geqslant 1$ and the variables $u=(u_1,\ldots ,u_n)$ let
\beq\label{tplusovi_trig}
T_{[n]}^{\pm} (u )=T_1^\pm( u_1)\ldots T_n^\pm ( u_n)
\fand
T_{[n]}^{\pm} (u|z)=T_1^\pm(z+u_1)\ldots T_n^\pm (z+u_n).
\eeq 
Finally, we recall  Etingof--Kazhdan's construction \cite[Thm. 2.3]{EK}.

\begin{thm}\label{EK:qva}
For any $c\in\CC,$ there exists a unique $h$-adic quantum vertex algebra  structure on $\Vc $
  such that 
	the vertex operator map $Y(\cdot ,z)$ is given by
\beq\label{qva1}
Y\big(T^+_{[n]}  (u)\vac,z\big)=T^+_{[n]}  (u|z)\ts T_{[n]}^- (u|z+hc/2)^{-1}, 
\eeq
the vacuum vector is $\vac $	
and the braiding  map $\mathcal{S}$ is defined by the relation  
\begin{align}
&\mathcal{S}(z)\big(R_{nm}^{  12}(e^{z+u-v})^{-1}  T_{[m]}^{+ 24}(v) 
R_{nm}^{  12}(e^{z+u-v-h  c})  T_{[n]}^{+ 13}(u)(\vac\otimes \vac) \big)\non\\
 =\, & T_{[n]}^{+ 13}(u)  R_{nm}^{  12}(e^{z+u-v+h  c})^{-1} 
 T_{[m]}^{+ 24}(v)  R_{nm}^{  12}(e^{z+u-v})(\vac\otimes \vac).\label{braiding}
\end{align}
\end{thm}

The next lemma gives a multiplicative counterpart of  \eqref{braiding}; cf.  \cite[Lemma 2.6]{K}.

\begin{lem}\label{shat_lemma_trig}
For any $c\in\CC ,$  there exists a unique $\CC[[h]]$-module map
$$
\wht{\Sc}(z)\colon \Vc\ot\Vc\to\Vc\ot\Vc\ot\CC(z)[[h]]
$$
which satisfies
\begin{align}
&\wht{\Sc}(z)\big(R_{nm}^{  12}(ze^{u-v})^{-1}  T_{[m]}^{+ 24}(v) 
R_{nm}^{  12}(ze^{u-v-h  c})  T_{[n]}^{+ 13}(u)(\vac\otimes \vac) \big)\non\\
 =\, & T_{[n]}^{+ 13}(u)  R_{nm}^{  12}(ze^{u-v+h  c})^{-1} 
 T_{[m]}^{+ 24}(v)  R_{nm}^{  12}(ze^{u-v})(\vac\otimes \vac).\label{Shatmaptrig}
\end{align}
\end{lem}

%%%%%%%%%%%%%%%%%%%%%%%%%%%%%%%%%%%
\subsection{Elliptic case}\label{sec0302}
%%%%%%%%%%%%%%%%%%%%%%%%%%%%%%%%%%%
In this subsection, we consider the elliptic $R$-matrix    \eqref{Ruh}.
As with the trigonometric case, we start by introducing  a certain algebra $\textrm{U}(R)=\textrm{U}(R^{ell}).$
In contrast with the papers \cite{EK3,EK4}, which motivated its definition,  it contains two distinct families of generators, which both, in the $h$-adic quantum vertex algebra construction below, play the role of creation operators.   
We  shall use the abbreviation 
$\eps_{12}=(-1)^{\eps_1-\eps_2} $ for $\eps_1,\eps_2\in\left\{0,1\right\}.$
Let $\ur$ be the topologically free associative algebra over the ring  $\CC[[h]]$ generated by the elements
$
t_{ij}^{(\eps,-r)},$ where $i,j=1,2,$ $\eps=0,1,$ $r=1,2,\ldots,$
subject to the defining relations 
\beq\label{frt}
R_{12}(\eps_{12}e^{u_1-u_2})\ts T_{1}^{\eps_1}(u_1)\ts T_{2}^{\eps_2}(u_2)=
T_{2}^{\eps_2}(u_2)\ts T_{1}^{\eps_1}(u_1)\ts R_{12}(\eps_{12}e^{u_1-u_2})
\quad\text{for}\quad
  \eps_1,\eps_2=0,1,
\eeq
where the matrices $T^\eps(u),$ $\eps=0,1,$  are given by
$$
T^\eps(u)=\sum_{i,j=1}^2 e_{ij}\ot t_{ij}^{\eps}(u)
\qquad\text{for} \qquad
t_{ij}^{\eps}(u)=\delta_{ij}-\sum_{r\geqslant 1}t_{ij}^{(\eps,-r)}\ts u^{r-1}.
$$

\begin{rem}\label{remdef_1}
As the matrices $T^\eps(u)$ contain only nonnegative powers of $u$,  the FRT-relations   \eqref{frt} are equivalent to
$$
R_{12}(\eps_{12}e^{-u_2+u_1})\ts T_{1}^{\eps_1}(u_1)\ts T_{2}^{\eps_2}(u_2)=
T_{2}^{\eps_2}(u_2)\ts T_{1}^{\eps_1}(u_1)\ts R_{12}(\eps_{12}e^{-u_2+u_1})
\quad\text{for}\quad
  \eps_1,\eps_2=0,1.
	$$
\end{rem}

\begin{rem}\label{remdef_2}
Note that \eqref{frt} gives only three families of defining relations, as the pairs $(\eps_1,\eps_2)=(0,1)$ and  $(\eps_1,\eps_2)=(1,0)$ yield equivalent relations. Indeed, this  is easily verified by employing the $R$-matrix property  \eqref{uniz}.
\end{rem}

\begin{rem}
 The rational counterpart of the defining relations \eqref{frt}, in the sense of Remark \ref{trig_ell_rem}, consists of two copies of the defining relations for the dual Yangian of $\mathfrak{gl}_2$ (for more details on the dual Yangian see, e.g., the paper by Iohara \cite{I}),
$$
R^{rat}_{+\ts 12} (u-v)\ts T_1^\eps(u)\ts T_2^\eps(v)=T_2^\eps(v)\ts T_1^\eps(u)\ts R^{rat}_{+\ts 12} (u-v)\quad\text{for}\quad \eps=0,1.
$$
 In addition, there are two equivalent families of relations
$$
R^{rat}_{-\ts 12} (u-v)\ts T_1^\eps(u)\ts T_2^{1-\eps}(v)=T_2^{1-\eps}(v)\ts T_1^\eps(u)\ts R^{rat}_{-\ts 12} (u-v)\quad\text{for}\quad \eps=0,1.
$$
\end{rem}

The following proposition is a simple consequence of Remarks \ref{remdef_1} and \ref{remdef_2}.

\begin{pro}\label{pro_23}
There exists a unique involutive automorphism $\tau$ of $\ur$ such that 
\beq\label{tau_formula}
\tau\colon t_{ij}^{(\eps,-r)}\mapsto t_{ij}^{(1-\eps,-r)}
\quad\text{for all}\quad i,j=1,2,\, \eps=0,1,\,r=1,2,\ldots .
\eeq
\end{pro}

Denote by $\vac$  the unit in   $\ur .$  From now on, we  consider   $T^\eps(u)$  as   operator series  over $\ur, $
such that its action is given by the algebra multiplication.  
We  proceed towards the construction of an $h$-adic quantum vertex algebra structure over the $\CC[[h]]$-module of  $\ur.$ First, we adapt the construction of annihilation operators from \cite[Lemma 2.1]{EK}.

\begin{lem}\label{lemma_21}
For any $c\in\CC ,$ there exists a unique invertible operator series
$$
T^-(u)\in\ndo\CC^2\ot\om(\ur,\ur((u))_h)
$$
such that for all $n\geqslant 0$ and $\eps=(\eps_1,\ldots ,\eps_n)$ we have
\begin{align}
&T^-_0(u)\ts T_{1}^{\eps_1}(v_1)\ldots T_{n}^{\eps_n}(v_n)\vac\non\\
=\,&
R^{01}_{1n}(e^{u-v^\eps+hc/2})^{-1}\ts T_{1}^{\eps_1}(v_1)\ldots T_{n}^{\eps_n}(v_n)\ts    R^{01}_{1n}(e^{u-v^\eps-hc/2})\vac.\label{lemma_21_eq}
\end{align}
\end{lem}

\begin{prf}
As with \cite[Lemma 2.1]{EK}, the lemma follows by a direct computation which relies on the defining relations \eqref{frt} and the Yang--Baxter equation \eqref{ell_addi}.
\end{prf}

From now on, we  denote $\ur$ by $\vr$ to indicate that it is equipped with the action of the operator series from Lemma \ref{lemma_21}, which depends on the   complex parameter $c$.
Denote by $T^-(u+\pi\iii)$ the operator series obtained from $T^-(u)$ by changing the signs of the arguments of  $R$-matrices in \eqref{lemma_21_eq}, so that  its action   is given by
\begin{align}
& T^-_0(u+\pi\iii)\ts T_{1\ts n+1}^{\eps_1}(v_1)\ldots T_{n\ts n+1}^{\eps_n}(v_n)\non\\
=\,&
R^{01}_{1n}(e^{u-v^{\eps^\prime}+hc/2})^{-1}\ts T_{1}^{\eps_1}(v_1)\ldots T_{n}^{\eps_n}(v_n)\ts R_{0n}^+\ldots  R^{01}_{1n}(e^{u-v^{\eps^\prime}-hc/2})\vac,\label{lemma_21_eq_v2}
\end{align}
where $ \eps^\prime  =(1-\eps_1,\ldots ,1-\eps_n).$
As with $T^-(u)$, this is   an invertible     series
in $\ndo\CC^2\ot\om(\ur,\ur((u))_h).$ Moreover, it exhibits the following property, which is an immediate consequence of \eqref{tau_formula}, \eqref{lemma_21_eq} and \eqref{lemma_21_eq_v2}.
\begin{lem}
We have
\beq\label{tauT}
T^-(u)^{\pm 1}\ts (1\ot\tau)=(1\ot\tau)\ts T^-(u+\pi\iii)^{\pm 1}.
\eeq
\end{lem}

The next lemma follows by a straightforward computation which relies on
 the defining relations \eqref{frt} and Lemma \ref{lemma_21}. 

\begin{lem}
The operator series $T^\eps(u)$ and $T^- (u)$ satisfy the FRT-relations
\begin{gather*}
R( e^{u_1-u_2})\ts T^-_{1} (u_1)\ts T^-_{2} (u_2)=
T^-_{2} (u_2)\ts T^-_{1} (u_1)\ts R( e^{u_1-u_2}),\\
R((-1)^\eps e^{u_1-u_2+hc/2})\ts T^-_{1} (u_1)\ts T_{2}^{\eps }(u_2)=
T_{2}^{\eps }(u_2)\ts T^-_{1} (u_1)\ts R((-1)^\eps e^{u_1-u_2-hc/2}).
\end{gather*}
\end{lem}

Finally, for $u$, $\eps$ and $u^\eps$ as in \eqref{uovi}--\eqref{uepsovi}, we shall use the following notation for the operator series with   coefficients in $(\ndo\CC^2)^{\ot n}\ot\ndo\vr$:
\begin{gather*}
 T_{[n]}^{\eps}(u) = T_{1}^{\eps_1}(u_1)\ldots T_{n}^{\eps_n}(u_n),\qquad T_{[n]}^{\eps}(u|z) = T_{1}^{\eps_1}(z+u_1)\ldots T_{n}^{\eps_n}(z+u_n),\\
 T^-_{[n]} (u^\eps |z )= T^-_{1} (z+u_1+ \eps_1\pi\iii )\ldots T^-_{n} (z+u_n+ \eps_n\pi\iii ) .
\end{gather*}

The next theorem slightly extends the elliptic case of the  Etingof--Kazhdan construction \cite[Thm. 2.3]{EK}.

\begin{thm}\label{thm_qva}
For any $c\in\CC ,$ there exists a unique structure of $h$-adic quantum vertex algebra over $\vr$ such that  the vertex operator map $Y(\cdot ,z)$ is given by
\beq\label{Ymap}
Y(T_{[n]}^{\eps}(u)\vac, z )
=T_{[n]}^{\eps}( u|z)\ts  T^-_{[n]} ( u^\eps|z +hc/2)^{-1} ,
\eeq
the vacuum vector is $\vac  $
and the braiding  map $\mathcal{S}$ is defined by the relation  
\begin{align}
&\Sc(z)\left(R_{nm}^{12}(e^{z+u^\eps-v^\nu})^{-1}\ts T_{[m]}^{\nu\ts 24}(v)\ts
R_{nm}^{12}(e^{z+u^\eps-v^\nu -hc})\ts T_{[n]}^{\eps\ts 13}(u)(\vac\ot\vac)
\right)\non\\
=\, &  T_{[n]}^{\eps\ts 13}(u)\ts R_{nm}^{12}(e^{z+u^\eps-v^\nu +hc})^{-1}\ts
T_{[m]}^{\nu\ts 24}(v)\ts R_{nm}^{12}(e^{z+u^\eps-v^\nu})(\vac\ot\vac).\label{Smap}
\end{align}
\end{thm}

\begin{prf}
We omit the   details, as they go in parallel with the proof of \cite[Thm. 2.3]{EK}. More specifically, one can  verify the $h$-adic quantum vertex algebra axioms from  \cite[Def. 2.20]{Li} by a direct calculation. In addition, all details can be also recovered by following the proofs of \cite[Thm. 2.3.8]{Gardini} and \cite[Thm. 4.1]{JKMY}.
\end{prf}

As with the trigonometric case (recall Lemma \ref{shat_lemma_trig}), we have the following simple lemma.

\begin{lem}\label{kor_shat_ell} 
For any $c\in\CC$  there exists a unique $\CC[[h]]$-module map
$$
\wht{\Sc}(z)\colon \vr\ot \vr\to \vr\ot \vr\ot \CC(z)[[h]] 
$$
which
satisfies
\begin{align}
&\wht{\Sc}(z)\left(R_{nm}^{12}(ze^{u^\eps-v^\nu})^{-1}\ts T_{[m]}^{\nu\ts 24}(v)\ts
R_{nm}^{12}(ze^{u^\eps-v^\nu -hc})\ts T_{[n]}^{\eps\ts 13}(u)(\vac\ot\vac)
\right)\non\\
=\, &  T_{[n]}^{\eps\ts 13}(u)\ts R_{nm}^{12}(ze^{u^\eps-v^\nu +hc})^{-1}\ts
T_{[m]}^{\nu\ts 24}(v)\ts R_{nm}^{12}(ze^{u^\eps-v^\nu})(\vac\ot\vac).\label{Shatmap}
\end{align}
\end{lem}

\begin{rem}
Recall the embedding \eqref{iotau_emb}.
The maps \eqref{Smap} and \eqref{Shatmap} are related by
\beq\label{SShat}
\Sc(z) = \iota_z\left(\wht{\Sc}(x)\big|_{x=e^z}\big.\right).
\eeq
The same applies to the trigonometric case, i.e., to the maps \eqref{braiding} and \eqref{Shatmaptrig}.
\end{rem}

The next definition is a straightforward generalization of the notion of vertex algebra homomorphism; see, e.g., the books by Kac \cite{Kac} and Lepowsky and Li \cite{LL}.

\begin{defn}\label{def_homo}
Let $(V^{(i)},Y^{(i)},\vac^{(i)},\Sc^{(i)})$, $i=1,2,$ be $h$-adic quantum vertex algebras. A $\CC[[h]]$-linear map $\psi\colon V^{(1)}\to V^{(2)}$ is said to be a {\em homomorphism} of $h$-adic quantum vertex algebras if it satisfies
$$
\psi\left(Y^{(1)}(v,z)w\right)=Y^{(2)}(\psi(v),z)\psi(w)\quad \text{for all}\quad v,w\in V^{(1)}\Fand
\phi(\vac^{(1)})=\vac^{(2)} .
$$
\end{defn}

\begin{pro}\label{propozicija316} 
The map $\tau$, as given by Proposition \ref{pro_23}, is an involutive automorphism of the $h$-adic quantum vertex algebra $\vr.$
\end{pro}

\begin{prf}
To verify the proposition, it suffices to show that the map $\tau$ is an $h$-adic vertex algebra homomorphism. Indeed, the remaining assertions are already evident from Proposition \ref{pro_23}. By Definition \ref{def_homo}, it is sufficient to prove that for any  families of variables $u$ and $v$ and   binary tuples $\eps$ and $\nu$, as in \eqref{uovi} and \eqref{epsovi}, we have
$$
\tau\left(
Y(T_{[n]}^{\eps\ts 13} (u)\vac,z)\ts T_{[m]}^{\nu\ts 23} (v)\vac
\right)
=
Y\left(\tau\left(T_{[n]}^{\eps\ts 13} (u)\vac\right) ,z\right)\ts \tau\hspace{-1pt}\left(T_{[m]}^{\nu\ts 23} (v)\vac
\right).
$$
By \eqref{Ymap}, the left-hand side equals
\beq\label{lhs}
\tau\left(
 T_{[n]}^{\eps\ts 13} (z+u) \ts  T_{[n]}^{-\ts 13} (z+u^\eps +hc/2)^{-1}\ts T_{[m]}^{\nu\ts 23} (v)\vac
\right).
\eeq
On the other hand, by using \eqref{tau_formula} and then \eqref{Ymap}, we find that the right-hand side  equals
\begin{align}
&Y( T_{[n]}^{ \eps^\prime\ts 13} (u)\vac  ,z)\ts   T_{[m]}^{ \nu^\prime\ts 23} (v)\vac 
= 
  T_{[n]}^{ \eps^\prime\ts 13} (z+u) \ts T_{[n]}^{-\ts  13} (z+u^{ \eps^\prime}+hc/2)^{-1} \ts  T_{[m]}^{ \nu^\prime\ts 23} (v)\vac.\label{rhs}
\end{align}
Finally, it follows from  \eqref{tau_formula} and \eqref{tauT}  that \eqref{lhs} and \eqref{rhs} coincide, as required.
\end{prf}

%%%%%%%%%%%%%%%%%%%%%%%%%%%%%%%%%%%
%%%%%%%%%%%%%%%%%%%%%%%%%%%%%%%%%%%
\section{Compatible pair \texorpdfstring{$(\sigma,\rho)$}{(sigma,rho)} associated with \texorpdfstring{$\vr$}{Vc(R)}}\label{sec04}
%%%%%%%%%%%%%%%%%%%%%%%%%%%%%%%%%%%
%%%%%%%%%%%%%%%%%%%%%%%%%%%%%%%%%%%

 In this section, we introduce a certain     pair of maps over the tensor square of 
$\vr=\mathcal{V}^{c}(R^{trig}),\mathcal{V}^{c}(R^{ell}).$
To
 consider both settings simultaneously, in the trigonometric case,  
we   denote the generator matrix $T^+(u)$,  defined by \eqref{tplustrig}, by $T^\eps(u).$ Furthermore, as  indicated in Subsection \ref{sec0103}, we assume that all binary tuples are trivial in the trigonometric case.
 First, for reader's convenience, we recall  
	a special case of the definition of multiplicative compatible pair \cite[Def. 3.6]{BK}, which fits the setting of this paper. 

\begin{defn}\label{def_br_mult}
Let $V$ be a topologically free $\CC[[h]]$-module and
$$
\sigma(z),\rho(z)\colon V\ot V\to V\ot V\ot \CC(z)[[h]] 
$$
  $\CC[[h]]$-module maps.
The pair $(\sigma, \rho)$ is said to be a {\em  multiplicative compatible pair} if it possesses the following properties.

\begin{enumerate}
\item The map $\sigma$ satisfies the {\em quantum Yang--Baxter equation},
\begin{align}
&\sigma_{12}(z_1  )\ts\sigma_{13}(z_1z_2    )\ts\sigma_{23}(z_2  )=
\sigma_{23}(z_2  )\ts\sigma_{13}(z_1z_2   )\ts\sigma_{12}(z_1 )  \label{ybe_muh}\\
\intertext{and  the {\em unitarity condition},}
&\sigma(1/z )\ts \sigma_{21}(z )=\sigma_{21}(z )\ts\sigma(1/z ) =1 .\label{uni_muh}
\end{align}

\item The map $\rho $
is invertible, i.e., there exists a map
$$   
\rho^{-1} (z )\colon V\ot V\to V\ot V\ot\CC    (z ) [[h]]
$$   
such that
$
\rho  (z )\ts\rho^{-1} (z )=\rho^{-1} (z )\ts \rho  (z )=1 .
$
\item The map 
\begin{align} 
 \mathcal{M}_{\sigma,\rho}  (z )\colon V\ot V&\to V\ot V\ot\CC   (z )[[ h]],\non\\
\mathcal{M}_{\sigma,\rho}  (z )&\coloneqq \rho (z )\ts \sigma (z )\ts \rho_{21}^{-1}(1/z )\label{map_m}
\end{align}
satisfies the    Yang--Baxter equation \eqref{ybe_muh} and the unitarity condition \eqref{uni_muh}.
\end{enumerate}
\end{defn}

  From now on, we omit the term ``multiplicative'' and refer to any pair of maps satisfying Definition \ref{def_br_mult} more briefly as a compatible pair.
In the next proposition, we give an example of  compatible pair over the $\CC[[h]]$-module of $\vr$. Its form resembles the one of the {\em additive} compatible pair given by \cite[Prop. 4.3]{BK0}. It is worth noting that the proposition makes use of the fact that the $R$-matrix $R(z)$ (both in the elliptic   and   trigonometric  case) can be   regarded as a formal rational function with respect to the variable $z$, i.e. as an element of $\ndo\CC^N\ot\ndo\CC^N(z)[[h]]$ (where $N=2$ in the elliptic case), which follows from  \cite[Rem. 2.4]{BBK} and \cite[Rem. 2.5]{K}.

\begin{pro}\label{lemma_31}
There exist unique $\CC[[h]]$-modules maps
$$
\sigma(z),\rho(z)\colon\vr\ot\vr\to\vr\ot\vr\ot\CC(z)[[h]]
$$
such that 
\begin{align}
&\sigma(z)\left(T_{[n]}^{\eps\ts 13}(u)\ts T_{[m]}^{\nu\ts 24}(v)(\vac\ot\vac)\right)
=
R_{nm}^{12}(ze^{u^\eps-v^\nu})\ts 
T_{[n]}^{\eps\ts 13}(u)\ts 
T_{[m]}^{\nu\ts 24}(v)\ts 
R_{nm}^{12}(ze^{u^\eps-v^\nu})^{-1}
(\vac\ot\vac),\label{sigma_formula}\\
&\rho(z)\left(T_{[n]}^{\eps\ts 13}(u)\ts T_{[m]}^{\nu\ts 24}(v)(\vac\ot\vac)\right)
=
T_{[n]}^{\eps\ts 13}(u)\ts 
R_{nm}^{12}(ze^{u^\eps-v^\nu+hc})^{-1}\ts
T_{[m]}^{\nu\ts 24}(v)\ts 
R_{nm}^{12}(ze^{u^\eps-v^\nu}) 
(\vac\ot\vac).\label{rho_formula}
\end{align}
Moreover, the maps $\sigma$ and $\rho$ form a compatible pair.
\end{pro}

\begin{prf}
The proposition can be proved by   directly verifying the constraints (1)--(3) imposed by Definition \ref{def_br_mult}. The fact  that the maps $\sigma$ and $\rho$ are well-defined by \eqref{sigma_formula} and \eqref{rho_formula} is established by arguing as in the proofs of \cite[Propositions 3.5, 3.9]{BK}.\vspace{1pt}

\noindent(1)
The Yang--Baxter equation \eqref{ybe_muh} for the map $\sigma$ is a consequence of the fact that the corresponding $R$-matrices \eqref{erofzed}  and \eqref{Rzh} satisfy the same identity; recall \eqref{trig_ybe_12} and \eqref{ell_ybe_12}.
As for the unitarity condition \eqref{uni_muh}, the trigonometric $R$-matrix \eqref{rbar_trigg} possesses the unitarity property $\R(z^{-1})=\R_{21}(z)^{-1},$ while the elliptic $R$-matrix \eqref{Rzh} satisfies the unitarity-like identity \eqref{uniz}. In particular, the $R$-matrices \eqref{erofzed}  and \eqref{Rzh}, which appear in the defining expression \eqref{sigma_formula} for $\sigma,$ are not unitary. However, the aforementioned properties still imply the  unitarity condition \eqref{uni_muh}, as the additional $R$-matrix terms cancel, due to the specific form of \eqref{sigma_formula}.
\vspace{1pt}

\noindent(2)
The invertibility of the map $\rho$ can be proved by writing the explicit formula for its inverse. For example, in the  trigonometric  case, we have
\beq\label{roalfabeta}
\rho(z)^{-1} = \beta(z)\ts \alpha(z) ,
\eeq
where
the maps
$
\alpha(z),\beta(z)\colon\vr\ot\vr\to\vr\ot\vr\ot\CC(z)[[h]]
$
are given by
\begin{align*}
&\alpha(z)\left(T_{[n]}^{\eps\ts 13}(u)\ts T_{[m]}^{\nu\ts 24}(v)(\vac\ot\vac)\right) \\
=&\,\,
T_{[n]}^{\eps\ts 13}(u)\ts
D^1_{[n]}\ts
R_{nm}^{12}(ze^{u^\eps-v^\nu+h(c+N)}) \ts
(D^1_{[n]})^{-1}\ts
 T_{[m]}^{\nu\ts 24}(v) 
(\vac\ot\vac),\\
&\beta(z)\left(T_{[n]}^{\eps\ts 13}(u)\ts T_{[m]}^{\nu\ts 24}(v)(\vac\ot\vac)\right) \\
=&\,\,
T_{[n]}^{\eps\ts 13}(u)\ts
 T_{[m]}^{\nu\ts 24}(v)\ts 
R_{nm}^{12}(ze^{u^\eps-v^\nu})^{-1} 
(\vac\ot\vac) 
\end{align*}
with $D^1_{[n]}=D^{\ot n} \ot 1^{\ot m}$; see \eqref{dijagonalna} for the definition of the matrix $D.$ To check that \eqref{roalfabeta} defines the inverse of $\rho,$  one needs to employ the first crossing symmetry property in  \eqref{csym}.
In the elliptic case, the explicit formula for the inverse of $\rho$ can be again derived by using the crossing symmetry property, as given by
\cite[Eq. 2.13]{FIJKMY2}.
\vspace{1pt}

\noindent(3) The fact that the map $\mathcal{M}_{\sigma,\rho},$  defined by \eqref{map_m},  
satisfies the    Yang--Baxter equation \eqref{ybe_muh} and the unitarity condition \eqref{uni_muh} is established in Remark \ref{remdef} below.
\end{prf}

Throughout the paper, we are interested in compatible pairs on $h$-adic quantum vertex algebras which are in tune with the underlying braiding map in the following sense.

\begin{defn}
Let $(V,Y,\Sc,\vac)$ be an $h$-adic quantum vertex algebra and $(\sigma,\rho)$ a compatible pair over $V$. The pair $(\sigma,\rho)$ is said to be {\em associated with the $h$-adic quantum vertex algebra $V$}  if the map $\mathcal{M}_{\sigma,\rho}$,   defined by \eqref{map_m}, satisfies 
\begin{align}\label{shtsm}
\Sc(z)=\iota_z \left(\mathcal{M}_{\sigma,\rho}(x)\big|_{x=e^z} \right).\big.
\end{align}
\end{defn}

Consider  the $h$-adic quantum vertex algebra $\Vc $ and the compatible pair from Proposition \ref{lemma_31}. To verify the requirement \eqref{shtsm}, it is sufficient to check that the map $\wht{\Sc},$ as given by Lemmas \ref{shat_lemma_trig} and \ref{kor_shat_ell}, satisfies
the identity 
\beq\label{s_hat_identity}
\wht{\Sc}(x)=  \mathcal{M}_{\sigma,\rho}(x). 
\eeq
This follows by a straightforward computation, so that we have the following corollary.

\begin{kor}\label{korolar44}
The pair $(\sigma,\rho),$ given by Proposition \ref{lemma_31}, is  associated with $\Vc.$
\end{kor}

\begin{rem}\label{remdef}
Note that, in  particular,  the identity 
\eqref{s_hat_identity} implies that the map $\mathcal{M}_{\sigma,\rho}$ satisfies 
the    Yang--Baxter equation \eqref{ybe_muh} and the unitarity equation \eqref{uni_muh}, as required by Definition \ref{def_br_mult}  as these equations hold for the braiding map $\wht{\Sc}.$
\end{rem}

Applying the embedding $\iota_z$ to the map $\sigma=\sigma(z),$ defined by \eqref{sigma_formula}, we obtain the map
\beq\label{iota-sigma}
\iota_z\sigma(z)\colon\vr\ot\vr\to\vr\ot\vr\ot\CC((z))[[h]].
\eeq
Next, applying the substitution $x=e^z$ to the map $\rho=\rho(x),$ defined by \eqref{rho_formula}, we get
$$
  \rho(x)\left|_{x=e^z} 
\colon\vr\ot\vr\to\vr\ot\vr\ot\CC_*(z)[[h]]
.\right.
$$
Denote by $\rho(e^z)$ its composition  with the embedding $\iota_z,$
\beq\label{rhoe}
\rho(e^z)\coloneqq \iota_z \left(\rho(x)\left|_{x=e^z}\right)
\colon\vr\ot\vr\to\vr\ot\vr\ot\CC((z))[[h]]
. \right.
\eeq
Observe that the composition  of the vertex operator map $Y$,  given by \eqref{qva1} and \eqref{Ymap}, and the map \eqref{rhoe} is well-defined. We shall denote it by $Y^\rho,$ so that we have
\beq\label{yrho_deff}
Y^\rho(z)\coloneqq Y(z)\ts\rho(e^z)\colon\vr\ot\vr\to\vr((z))_h\ot\CC((z))[[h]].
\eeq
 The above maps \eqref{iota-sigma} and  \eqref{rhoe} satisfy the following hexagon-type identities.
\begin{pro}\label{propozicija45}
We have
\begin{align}
&\iota_{z_1}\sigma( z_1 )\left(Y^\rho(z_2)\ot 1\right)=
\left(Y^\rho(z_2)\ot 1\right) \iota_{z_1}\sigma_{23}( z_1 )  \left(\iota_{z}\sigma_{13}(z)\right)\big|_{z=z_1 e^{z_2}},\big. \label{hex_sigma}\\
&\rho(e^{z_1})\left(Y^\rho(z_2)\ot 1\right)=
\left(Y^\rho(z_2)\ot 1\right) \rho_{13}(e^{z_1+z_2})\ts\rho_{23}(e^{z_1}).\label{hex_rho}
\end{align}
\end{pro}

\begin{prf}
Clearly, the coefficients of the matrix entries of all expressions of the form 
\beq\label{expressionn}
T_{[n]}^{\eps\ts 14}(u)\ts 
T_{[m]}^{\nu\ts 25}(v)\ts 
T_{[k]}^{\mu\ts 36}(w) 
\eeq
  form an $h$-adically dense $\CC[[h]]$-submodule of $\vr^{\ot 3}$.
Hence, to prove the proposition, it suffices to show that the images of   \eqref{expressionn} under the both sides of \eqref{hex_sigma}  and \eqref{hex_rho} coincide. 

First, by employing   the definitions of the   maps $Y$ and $\rho,$   given by  
\eqref{qva1}, \eqref{Ymap} and  \eqref{rhoe},
 one obtains the explicit formula for the action of      \eqref{yrho_deff},  
\beq\label{yrho_formula}
Y^\rho (  z )\, T_{[n]}^{\eps\ts 13}(u)\ts 
T_{[m]}^{\nu\ts 24}(v)
=
T_{[n]}^{\eps\ts 13}(u|z)\ts 
T_{[m]}^{\nu\ts 23}(v).
\eeq
By using the definition 
\eqref{sigma_formula} of the map $\sigma$ and \eqref{yrho_formula},
one   finds that the images of \eqref{expressionn} under the both sides of the identity \eqref{hex_sigma} are equal to
$$
R_{n+m\ts k}^{12\, 3}(z_1 e^{x^{\eta} -w^\mu})\ts
T_{[n]}^{\eps\ts 14}(u|z_2)\ts 
T_{[m]}^{\nu\ts 24}(v)\ts 
T_{[k]}^{\mu\ts 35}(w)\ts 
R_{n+m\ts k}^{12\, 3}(z_1 e^{x^{\eta} -w^\mu})^{-1},
$$
where
\beq\label{xandeta}
x=( z_2+u_1,\ldots , z_2+u_n ,v_1,\ldots ,v_m)\fand
\eta=(\eps,\nu)=(\eps_1,\ldots ,\eps_n,\nu_1,\ldots ,\nu_m).
\eeq
Analogously, by using the definition 
\eqref{rho_formula} of the map $\rho$ and \eqref{yrho_formula},
one   finds that the images of \eqref{expressionn} under the both sides of the identity \eqref{hex_rho}   equal  
$$
T_{[n]}^{\eps\ts 14}(u|z_2)\ts 
T_{[m]}^{\nu\ts 24}(v)\ts 
R_{n+m\ts k}^{12\, 3}( e^{z_1+x^{\eta} -w^\mu +hc})^{-1}\ts
T_{[k]}^{\mu\ts 35}(w)\ts 
R_{n+m\ts k}^{12\, 3}(  e^{z_1+x^{\eta} -w^\mu}) ,
$$
where $x$ and $\eta$ are  again given by \eqref{xandeta}. Thus, we conclude that both assertions of the proposition hold.
\end{prf}

In contrast with the rest of this section, the following simple corollary applies only to the elliptic case, where the (extended) $h$-adic quantum vertex algebra $\Vc$ admits the involutive automorphism $\tau$; recall Proposition \ref{propozicija316}. 
It  follows by comparing the defining expressions \eqref{tau_formula}, \eqref{sigma_formula} and \eqref{rho_formula} for the maps $\tau,$ $\sigma$ and $\rho.$  
\begin{kor}
For $i=1$ (resp. $i=2$) let $\tau_i$ denote the action of the automorphism $\tau$ on the first (resp. second) tensor factor of $\Vc\ot\Vc.$ We have 
$$
\tau_i\circ \sigma(z) =\sigma(-z)\circ \tau_i\fand
\tau_i\circ \rho(z) =\rho(-z)\circ \tau_i,
$$
where $i=1,2,$
and, consequently,
$$
\left(\tau\ot\tau\right) \circ \sigma(z) =\sigma( z)\circ \left(\tau\ot\tau\right)\fand
\left(\tau\ot\tau\right) \circ \rho(z) =\rho( z)\circ \left(\tau\ot\tau\right).
$$
\end{kor}

%%%%%%%%%%%%%%%%%%%%%%%%%%%%%%%%%%%
%%%%%%%%%%%%%%%%%%%%%%%%%%%%%%%%%%%
\section{Constructing deformed \texorpdfstring{$\phi$}{phi}-coordinated  \texorpdfstring{$\vr$}{Vcr}-modules}\label{sec05}
%%%%%%%%%%%%%%%%%%%%%%%%%%%%%%%%%%%
%%%%%%%%%%%%%%%%%%%%%%%%%%%%%%%%%%%

In this section, we employ   compatible pairs
to modify the notion of (un)twisted $\phi$-coordinated module for $h$-adic quantum vertex algebra. Finally, we investigate the connection between such structures and representation theories of the $h$-Yangian $\Y$ and the elliptic quantum algebra $\ahp.$  

%%%%%%%%%%%%%%%%%%%%%%%%%%%%%%%%%%%
\subsection{Modules over the \texorpdfstring{$h$}{h}-Yangian \texorpdfstring{$\Y$}{Yh(glN)}}\label{sec0501}
%%%%%%%%%%%%%%%%%%%%%%%%%%%%%%%%%%%
 Let
 $\phi=\phi(z_2,z_0)=z_2 e^{ z_0}\in\CC[[z_0,z_2]]$  be the {\em associate of the one-dimensional additive formal group}; see \cite[Sect. 2]{Liphi} for more details.
The next definition  combines the notions of {\em  $\phi$-coordinated module}   \cite[Def. 3.4]{Liphi} 
and {\em $(\sigma,\rho)$-deformed module} \cite[Def. 3.5]{BK0}.

\begin{defn}\label{def_mod_phi_trig}
Let $(V,Y,\vac,\Sc)$  be an $h$-adic    quantum vertex algebra and  $(\sigma   ,\rho   )$      a
   compatible pair associated with $V$.  
 Let  $W$  be a topologically free $\CC[[h]]$-module equipped with  a $\CC[[h]]$-module map
\begin{align*} 
Y_W(\cdot ,z) \colon V\ot W&\to W((z ))_h, \\
v\ot w&\mapsto Y_W(z)(v\ot w)=Y_W(v,z)w=\sum_{r\in\ZZ}  v_{r-1} w \ts z^{-r }.\non
\end{align*}
A pair $(W,Y_W)$ is said to be a {\em   $(\sigma   ,\rho   )$-deformed $\phi$-coordinated $V$-module} if the   map $Y_W(\cdot, z)$ satisfies 
\begin{align}
&Y_W(\vac,z)w=w\quad\text{for all }w\in W,\non
\intertext{the  {\em weak  $\rho  $-associativity}:
for any elements $u,v \in V$  and $n \in\mathbb{Z}_{> 0}$
there exists  $r\in\mathbb{Z}_{\geqslant 0}$
such that we have}
&(z_1-z_2)^r\ts Y_W(u,z_1)\ts Y_W(v,z_2)\in \om\left(W,W((z_1,z_2))\right) \mod h^n,\label{weakassoc_1}\\
&\left((z_1-z_2)^r\ts\ts Y_W(u, z_1)  Y_W(v,z_2)  \right)\Big|_{z_1=z_2 e^{z_0 }}^{\text{mod } h^n} \Big.\label{weakassoc_2}\\
&\qquad\big. -  z_2^r (e^{z_0} -1)^r\ts Y_W\big(Y^{\rho  }  (u, z_0  )v ,z_2\big)\,\in\, h^n \om\left(W,W[[z_0^{\pm 1},z_2^{\pm 1}]]\right)   ,\non
\intertext{and the  {\em $\sigma  $-locality}:
for any $u,v\in V$ and $n\in\mathbb{Z}_{> 0}$ there exists
    $r\in\mathbb{Z}_{\geqslant 0}$ such that for all $w\in W$ we have }
 &\big((z_1-z_2)^r\ts Y_W(z_1)\big(1\otimes Y_W(z_2)\big)\big( \sigma  (z_1 /z_2 )(u\otimes v)\otimes w\big)\big. 
\non\\
 &\qquad\big. - (z_1-z_2)^r\ts   Y_W(v,z_2) Y_W(u, z_1)  w \big)  \,
\in\,  h^n W[[z_1^{\pm 1},z_2^{\pm 1}]]. \label{sigma-loc-trigg}
\end{align}
\end{defn} 

\begin{rem}\label{remark_52}
In this remark, we explain some subtleties of  Definition \ref{def_mod_phi_trig}. First of all, by \eqref{weakassoc_1}, we require that the given expression, when applied on an arbitrary element $w\in W$ and then regarded modulo $h^n,$ possesses finitely many negative powers of $z_1$ and $z_2.$ In other words, we require that for any $w\in W$ we have
\beq\label{weakassoc_6}
(z_1-z_2)^r\ts Y_W(u,z_1)\ts Y_W(v,z_2)w = A_w(z_1,z_2) +  h^n B_w(z_1,z_2) 
\eeq
 for some formal series
$A_w(z_1,z_2)\in W((z_1,z_2))$ and $B_w(z_1,z_2)\in W[[z_1^{\pm 1},z_2^{\pm 1}]]$ (which both depend on $w$).
Note that this still does not ensure that one can apply the substitution $z_1=z_2 e^{z_0 }$ to \eqref{weakassoc_6}, as the resulting expression does not need to be well-defined. However, this substitution can be applied to $A_w(z_1,z_2).$ The notation used in  \eqref{weakassoc_2} indicates that the given expression is to be regarded modulo $h^n$ before the substitution is applied. In other words, we have
$$
\left((z_1-z_2)^r\ts\ts Y_W(u, z_1)  Y_W(v,z_2)  \right)\Big|_{z_1=z_2 e^{z_0 }}^{\text{mod } h^n}w=A_w(z_1,z_2).
$$
Finally, we remark that the form of the weak  $\rho  $-associativity property from Definition \ref{def_mod_phi_trig} is motivated by \cite[Rem. 3.2]{Liphi}; see also \cite[Lemma 2.9]{LiTW0}.
Regarding the $\sigma$-locality property \eqref{sigma-loc-trigg},
as with \eqref{iota-sigma},
   we assume that    the image of $\sigma$ is embedded in  $V\ot V\ot\CC((z_1  /z_2 )),$ i.e., we have $\sigma  (z_1 /z_2 ) =\left(\iota_z\ts\sigma  (z)\right)\big|_{z=z_1 /z_2 }.
$
\end{rem}

To associate  $\phi$-coordinated $(\sigma,\rho)$-deformed  $\vr$-modules with the $h$-Yangian, we introduce the notion of restricted   $\Y$-module. It is motivated by restricted modules for the affine Kac--Moody  Lie algebras; see, e.g., the book by Kac \cite{Kac_Inf}. A $\Y$-module $W$ is said to be  
{\em restricted} if it is topologically free as a $\CC[[h]]$-module and the corresponding action $L^-(z)_W$ of the generator matrix \eqref{elminusihY} is such that we have
$$
L^-(z)_W\in\ndo\CC^N \ot\om (W,W[z^{-1}]_h).
$$

  From now on, we denote by    $\sigma$ and $\rho$   the maps     \eqref{sigma_formula} and \eqref{rho_formula}, respectively.  

\begin{thm}\label{thm_53-trigg}
Let $W$ be a restricted $\Y$-module.
For any $c\in\CC,$ there exists a unique structure of $(\sigma   ,\rho    )$-deformed $\phi$-coordinated $\vr$-module over $W$ such that the module map is given by
\beq\label{thm53eq}
Y_{W}(T_{[n]}^{+}(u )\vac,z)
=L_{1}^- (z e^{u_1})_W \ldots L_{n}^- (ze^{u_n})_W .
\eeq 
Conversely, if $(W,Y_W)$ is a $(\sigma   ,\rho    )$-deformed $\phi$-coordinated $\vr$-module such that
\beq\label{assignment2}
Y_W(t_{ij}^+(0)\vac,z )\in \delta_{ij}+h\om(W,W[z^{-1}]_h)\quad\text{for all }i,j=1,\ldots ,N,
\eeq
then the assignment
\beq\label{assignment1}
L^-(z)_W\coloneqq Y_W(T^+(0)\vac,z )
\eeq
defines a structure of restricted $\Y$-module over W.
\end{thm}

\begin{prf}
Let $W$ be a restricted $\Y$-module. 
The fact that \eqref{thm53eq} defines a structure of $(\sigma   ,\rho    )$-deformed $\phi$-coordinated $\vr$-module over $W$ can be proved by straightforward computations which directly  verify the constraints imposed by Definition \ref{thm_53-trigg}.
We omit the details as they go in parallel with the proof of \cite[Thm. 5.4]{BK}, which is slightly more technical due to the form of the underlying compatible pair.

Conversely, suppose
$(W,Y_W)$ is a $(\sigma   ,\rho    )$-deformed $\phi$-coordinated $\vr$-module such that     \eqref{assignment2} holds.
To prove the second assertion of the theorem, it suffices to check that \eqref{assignment1} satisfies the defining relation \eqref{rll_def_trig} for $\Y$.
Using the   expression \eqref{sigma_formula} for the map $\sigma $
and the $\sigma$-locality \eqref{sigma-loc-trigg}, one finds that  for any $n\geqslant 1$ there exist $r\geqslant 0$ such that
$$
(z_1-z_2)^r R(z_1/z_2)\ts L^-_{1}(z_1)_W\ts L^-_{2}(z_2)_W\ts w \ts R(z_1/z_2)^{-1}
=(z_1-z_2)^r  L^-_{2}(z_2)_W\ts L^-_{1}(z_1)_W\ts w\mod h^n,
$$
where the map $L^-(z)_W$ is defined by \eqref{assignment1}.
Due to \eqref{assignment2}, the terms $(z_1 -z_2)^r$ cancel, so that  multiplying the above identity from the right by $R(z_1/z_2)$ yields
$$
  R(z_1/z_2)\ts L^-_{1}(z_1)_W\ts L^-_{2}(z_2)_W\ts w  
=   L^-_{2}(z_2)_W\ts L^-_{1}(z_1)_W \ts R(z_1/z_2)\ts w\mod h^n .
$$
Finally, with $W$ being separated, this implies 
$$
  R(z_1/z_2)\ts L^-_{1}(z_1)_W\ts L^-_{2}(z_2)_W\ts w  
=   L^-_{2}(z_2)_W\ts L^-_{1}(z_1)_W \ts R(z_1/z_2)\ts w 
$$
for all $w\in W ,$
as required.
\end{prf}

%%%%%%%%%%%%%%%%%%%%%%%%%%%%%%%%%%%
\subsection{Modules over the elliptic quantum algebra \texorpdfstring{$\ahp$}{Ahp(gl2)}}\label{sec0502}
%%%%%%%%%%%%%%%%%%%%%%%%%%%%%%%%%%%

In this subsection, we consider the associate
 $\phi=\phi(z_2,z_0)=z_2 e^{2 z_0}\in\CC[[z_0,z_2]]$.
Let $\omega_M=e^{2\pi\iii/M}$ be the principal primitive $M$-th root of unity.
The next definition  combines the notions of {\em twisted $\phi$-coordinated module}   \cite[Def. 2.8]{LiTW} 
and {\em $(\sigma,\rho)$-deformed module} \cite[Def. 3.5]{BK0}.

\begin{defn}\label{def_mod_phi}
Let $(V,Y,\vac,\Sc)$  be an $h$-adic    quantum vertex algebra, $\tau$  an automorphism of $V$ of period $M$ and $(\sigma   ,\rho   )$     a
   compatible pair associated with $V.$ 
Let  $W$  be a topologically free $\CC[[h]]$-module equipped with  a $\CC[[h]]$-module map
\begin{align*} 
Y_W(\cdot ,z) \colon V\ot W&\to W((z^{1/M}))_h, \\
v\ot w&\mapsto Y_W(z)(v\ot w)=Y_W(v,z)w=\sum_{r\in\frac{1}{M}\ZZ}  v_{r-1} w \ts z^{-r }.\non
\end{align*}
A pair $(W,Y_W)$ is said to be a {\em   $(\sigma   ,\rho ,\tau  )$-deformed $\phi$-coordinated $V$-module} if the   map $Y_W(\cdot, z)$ satisfies 
\begin{align}
&Y_W(\vac,z)w=w\quad\text{for all }w\in W,\non\\
&Y_W(\tau\ts u, z) =\lim_{x^{1/M}\to \omega_M^{-1}z^{1/M}}Y_W( u, x)\quad\text{for all }u\in V,\label{propertywtihtau}\\
\intertext{the  {\em weak  $\rho  $-associativity}:
for any elements $u,v \in V$  and $n \in\mathbb{Z}_{> 0}$
there exists  $r\in\mathbb{Z}_{\geqslant 0}$
such that we have}
&(z_1-z_2)^r\ts Y_W(u,z_1)\ts Y_W(v,z_2)\in \om\left(W,W((z_1^{1/M},z_2^{1/M}))\right) \mod h^n,\label{weakassoc_1_ellipt}\\
&\left((z_1-z_2)^r\ts\ts Y_W(u, z_1)  Y_W(v,z_2)  \right)\Big|_{z_1^{1/M}=(z_2 e^{2z_0 })^{1/M}}^{\text{mod } h^n} \Big.\non\\
&\qquad\big. -  z_2^r (e^{2z_0} -1)^r\ts Y_W\big(Y^{\rho  }  (u, z_0  )v ,z_2\big)\,\in\, h^n \om\left(W,W[[z_0^{\pm 1/M},z_2^{\pm 1/M}]]\right)   ,\label{weakassoc_2_ellipt}
\intertext{and the  {\em $\sigma  $-locality}:
for any $u,v\in V$ and $n\in\mathbb{Z}_{> 0}$ there exists
    $r\in\mathbb{Z}_{\geqslant 0}$ such that for all $w\in W$ we have }
 &\big((z_1-z_2)^r\ts Y_W(z_1)\big(1\otimes Y_W(z_2)\big)\big( \sigma  (z_1^{1/M} /z_2^{1/M} )(u\otimes v)\otimes w\big)\big. 
\non\\
 &\qquad\big. - (z_1-z_2)^r\ts   Y_W(v,z_2) Y_W(u, z_1)  w \big)  \,
\in\,  h^n W[[z_1^{\pm 1/M},z_2^{\pm 1/M}]]. \label{sigma-loc}
\end{align}
\end{defn}

In Definition \ref{def_mod_phi},
we assume that all technicalities explained in 
Remark \ref{remark_52}
suitably apply to the 
weak  $\rho  $-associativity property \eqref{weakassoc_1_ellipt}--\eqref{weakassoc_2_ellipt} and $\sigma  $-locality property \eqref{sigma-loc}.
In particular, as with \eqref{iota-sigma}, the image of $\sigma$ in \eqref{sigma-loc}  is supposed to be   embedded in  $V\ot V\ot\CC((z_1^{1/M} /z_2^{1/M})),$ i.e., we have
$$
\sigma  (z_1^{1/M} /z_2^{1/M} ) =\left(\iota_z\ts\sigma  (z)\right)\big|_{z=(z_1 /z_2)^{1/M} }.
$$

\begin{rem}
Setting $\sigma=\wht{\Sc}$ and $\rho=1 $ in Definition \ref{def_mod_phi}, one obtains   the $h$-adic version of the notion of   twisted $\phi$-coordinated module of Li, Tan and Wang   \cite[Def. 2.8]{LiTW}. 
\end{rem}

Our next goal is to construct a structure of $(\sigma   ,\rho ,\tau  )$-deformed $\phi$-coordinated $\vr$-module over $\ahp .$ From now on, we denote by $\tau,$  $\sigma$ and $\rho$   the maps given by \eqref{tau_formula}, \eqref{sigma_formula} and \eqref{rho_formula}, respectively. By  Proposition \ref{pro_23}, the automorphism $\tau$ is of period $2$, so   we consider Definition \ref{def_mod_phi} with $M=2.$ 

\begin{thm}\label{thm_53}
For any $c\in\CC,$ there exists a unique structure of $(\sigma   ,\rho ,\tau  )$-deformed $\phi$-coordinated $\vr$-module over $\ahp  $ such that the module map is given by
$$
Y_{\ahp}(T_{1}^{\eps_1}(u_1)\ldots T_{n}^{\eps_n}(u_n)\vac,z)
=L_{1} (z^{1/2}(-1)^{\eps_1}e^{u_1})\ldots L_{n} (z^{1/2}(-1)^{\eps_n}e^{u_n}).
$$ 
Consequently, if $W$ is an $\ahp  $-module which is topologically free as a $\CC[[h]]$-module, it is naturally equipped with a structure of $(\sigma   ,\rho ,\tau  )$-deformed $\phi$-coordinated $\vr$-module such that the module map is given by
\beq\label{YW_mod_formulla}
Y_{W}(T_{1}^{\eps_1}(u_1)\ldots T_{n}^{\eps_n}(u_n)\vac,z)
=L_{1} (z^{1/2}(-1)^{\eps_1}e^{u_1})_W \ldots L_{n} (z^{1/2}(-1)^{\eps_n}e^{u_n})_W .
\eeq
\end{thm}

\begin{prf}
As with Theorem
\ref{thm_53-trigg},
this theorem can be verified by direct arguments which closely follow the proof of \cite[Thm. 5.4]{BK}. Therefore, we only
prove the property
\eqref{propertywtihtau}, which involves the automorphism $\tau,$ as it does not appear in the aforementioned proof.
It is a simple consequence of
\eqref{tau_formula} and \eqref{YW_mod_formulla}. Indeed, if $W$ is an $\ahp  $-module, which is topologically free as a $\CC[[h]]$-module, we have
\begin{align*}
&\, Y_{W}(\tau (T_{1}^{\eps_1}(u_1)\ldots T_{n}^{\eps_n}(u_n)\vac),z)\\
=&\,Y_{W}( T_{1}^{1-\eps_1}(u_1)\ldots T_{n}^{1-\eps_n}(u_n)\vac ,z)\\
=&\,L_{1} (z^{1/2}(-1)^{1-\eps_1}e^{u_1})_W \ldots L_{n} (z^{1/2}(-1)^{1-\eps_n}e^{u_n})_W \\
=&\,L_{1} ((-1)z^{1/2}(-1)^{ \eps_1}e^{u_1})_W \ldots L_{n} ((-1)z^{1/2}(-1)^{ \eps_n}e^{u_n})_W\\
=&\,\lim_{x^{1/2}\to (-1) z^{1/2}}L_{1} (x^{1/2}(-1)^{ \eps_1}e^{u_1})_W \ldots L_{n} (x^{1/2}(-1)^{ \eps_n}e^{u_n})_W\\
=&\,\lim_{x^{1/2}\to (-1) z^{1/2}} Y_{W}( T_{1}^{\eps_1}(u_1)\ldots T_{n}^{\eps_n}(u_n)\vac ,x),
\end{align*} 
which   implies  \eqref{propertywtihtau},
as required.
\end{prf}

Suppose $(W,Y_W)$ is a $(\sigma   ,\rho ,\tau  )$-deformed $\phi$-coordinated $\vr$-module.
The   converse of the second assertion of Theorem \ref{thm_53} does not need to hold, due to the special form of the generator matrix $L(z)$; cf. \eqref{Lijn} and \eqref{Lijnn}. However, a partial converse can be established by adding suitable  requirements on the  map $Y_W(\cdot ,z).$ More specifically,  write
$$
Y_W(t_{ij}^{(0,-1)}\vac,z^2)=\sum_{n\in\ZZ} (t_{ij}^{(0,-1)}\vac )_n \ts  z^{-n},\quad\text{where } \quad i,j=1,2  .
$$
The coefficients $(t_{ij}^{(0,-1)}\vac )_n  $ of the variable $z$   belong to $ \ndo W $ for all $n.$ Suppose that their images satisfy 
\beq\label{condition_1}
\im \ts(t_{ij}^{(0,-1)}\vac )_n \subseteq h^{bn}W\quad\text{for all}\quad r> 0,
\eeq
where $b$ is a positive integer given by \eqref{pandh}.
This implies that the matrix entries of 
$$
\Lc(  z)_W \coloneqq Y_W(T^{(0)}(0)\vac,z^2)\in\ndo\CC^2\ot \om (W,W((z))_h)
$$
 satisfy \eqref{Lijn}. Furthermore, assume that
\beq\label{condition_2}
\ts(t_{ij}^{(0,-1)}\vac )_n =0\quad\text{if}\quad (-1)^{i+j}\neq (-1)^n.
\eeq
Then the matrix entries of the operator series $\Lc(  z)_W$ satisfy the condition  \eqref{Lijnn} as well. Finally, it remains to observe that  the $\sigma$-locality \eqref{sigma-loc}   implies the relation
$$
R(z_1/z_2)\ts \Lc_1(z_1)_W\ts \Lc_2(z_2)_W
=\Lc_2(z_2)_W\ts \Lc_1(z_1)_W\ts R(z_1/z_2).
$$
In particular, the requirement \eqref{condition_1} ensures that the matrix entries of $\Lc_1(z_1)_W  \Lc_2(z_2)_W$ and $\Lc_2(z_2)_W  \Lc_1(z_1)_W$ belong to $\om(W,W((z_1,z_2))_h),$
so that the 
terms $(z_1-z_2)^r$ in the $\sigma$-locality cancel. Finally, the preceding discussion implies the following theorem.

\begin{thm}\label{thm_66}
Let $(W,Y_W)$ be a $(\sigma   ,\rho ,\tau  )$-deformed $\phi$-coordinated $\vr$-module such that the map $Y_W(\cdot ,z)$ satisfies \eqref{condition_1} and  \eqref{condition_2}. Then the assignment 
$$
L(  z)_W \coloneqq Y_W(T^{0}(0)\vac,z^2)
$$
defines a structure of $\ahp$-module over $W.$  
\end{thm}

%%%%%%%%%%%%%%%%%%%%%%%%%%%%%%%%%%%
%%%%%%%%%%%%%%%%%%%%%%%%%%%%%%%%%%%
\section{Central elements in \texorpdfstring{$\mathcal{V}^{crit}(R^{trig})$}{Vcrit(Rtrig)}  and commutative families in   \texorpdfstring{$\Y$}{Yh(glN)}}\label{sec06}
%%%%%%%%%%%%%%%%%%%%%%%%%%%%%%%%%%%
%%%%%%%%%%%%%%%%%%%%%%%%%%%%%%%%%%%

In this section, we consider only the trigonometric case. We denote by $\mathcal{V}^{crit}(R)= \mathcal{V}^{crit}(R^{trig})$ the  $h$-adic quantum vertex algebra $\vr$ from Theorem \ref{EK:qva} associated with the trigonometric $R$-matrix \eqref{erofu} at the critical level $c=-N.$

 %%%%%%%%%%%%%%%%%%%%%%%%%%%%%%%%%%%
\subsection{Fusion procedure and fixed points of the   pair \texorpdfstring{$(\sigma,\rho)$}{(sigma,rho)}}\label{sec0601}
%%%%%%%%%%%%%%%%%%%%%%%%%%%%%%%%%%%

In this subsection, following the exposition in \cite[Sect. 2, 3]{JLM}, we recall some consequences of the fusion procedure for the Hecke
algebra \cite{C,IMO,N},   given in terms of the   $R$-matrices
$$
\Rv(z)=\frac{e^{h/2}-e^{-h/2}z}{1-z}\R(z)
\Fand
\Rv(e^u)=\frac{e^{h/2}-e^{u-h/2}}{1-e^u}\R(e^u),
$$
where $\R(z)$ is defined by \eqref{rbar_trigg} and $\R(e^u)=\left(\R(x)\right)\big|_{x=e^u}.$ Finally, we use them  to construct families of fixed points of the compatible pair $(\sigma, \rho) $ given by Proposition \ref{lemma_31}.

Suppose $\Lambda$ is a standard tableau of shape $\lambda\vdash n$ and denote by $c_k(\Lambda)$ the content $j-i$ of the box $(i,j)$ of $\lambda$ occupied by $k$ in $\Lambda.$
Consider the order over the set of all pairs $(i,j)$, where 
$1\leqslant i <j\leqslant n,$ given by
\beq\label{order}
(i,j )\prec (i^\prime, j^\prime)
\qquad\text{if}\qquad
j<j^\prime
\quad\text{or}\quad
j=j^\prime\text{ and }i<i^\prime .
\eeq
Let
$$
\Rv_\Lambda(z_1,\ldots ,z_n)
=
\prod_{(i,j) }^{\longrightarrow} 
\left(
P_{j-i\ts j-i+1}\ts \Rv_{j-i\ts j-i+1}(z_i e^{ (c_i(\Lambda)-c_j(\Lambda))h}/z_j)
\right),
$$
where the arrow indicates that the product over the set of all pairs $(i,j)$ with $1\leqslant i <j\leqslant n$ is ordered with respect to \eqref{order} and $P_{rs}$ is the action of the permutation operator $P=\sum_{i,j=1}^N   e_{ij}  \ot e_{ji}$
over the $r$-th and $s$-th tensor factor of $(\ndo\CC^N)^{\ot n}.$
Let $\lambda^\prime$ be the conjugate partition of $\lambda,$ $c_{\lambda^\prime}$ the Schur element associated with $\lambda^\prime$ and $\check{R}_0$ a certain invertible operator over $(\CC^N)^{\ot n}$; see\cite{JLM} for the details. By the fusion procedure,  
$$
\Ec_\Lambda = \frac{1}{c_{\lambda^\prime}}\ts \Rv_\Lambda(z_1,\ldots ,z_n)\ts
\check{R}_0^{-1}\Big|_{z_1=1}\Big|_{z_2=1}\ldots \Big|_{z_n=1},
$$ 
is a well-defined operator satisfying $\Ec_\Lambda^2=\Ec_\Lambda.$ 
 Clearly, the above identity  can be    written in terms of the $R$-matrices $\check{R}(e^u)$ as
	\begin{align}
&\Ec_\Lambda = \frac{1}{c_{\lambda^\prime}}\ts \Rv_\Lambda(e^{u_1},\ldots ,e^{u_n})\ts
\check{R}_0^{-1}\Big|_{u_1=0}\Big|_{u_2=0}\ldots \Big|_{u_n=0},\qquad\text{where}\label{fusion2}\\
&\Rv_\Lambda(e^{u_1},\ldots ,e^{u_n})
=
\Rv_\Lambda(z_1,\ldots ,z_n)
\Big|_{z_1=e^{u_1}}\Big|_{z_2=e^{u_2}}\ldots \Big|_{z_n=e^{u_n}}.\non
\end{align}
By using this observation, one easily derives  
the following properties of the idempotent $\Ec_\Lambda $ (see \cite[Lemmas 3.2, 3.3]{JLM}).

\begin{lem}\label{lemma_61}
Consider the $n$-tuples of variables   
\begin{align*}
&u^{(\Lambda)}=(u -c_1(\Lambda)h,\ldots ,u -c_n(\Lambda)h),\\
&u^{(\Lambda)}-Nh/2 =(u -c_1(\Lambda)h-Nh/2,\ldots ,u -c_n(\Lambda)h-Nh/2),\\
&z^{(\Lambda)}=(z  e^{-c_1(\Lambda)h},\ldots ,z  e^{-c_n(\Lambda)h}).
\end{align*}
\begin{enumerate}[(a)]
\item The following relations hold for operators on $\mathcal{V}^{crit}(R)$:
\begin{align}
&T_{[n]}^+ (u^{(\Lambda)} )  \ts \Ec_\Lambda =\Ec_\Lambda\ts  T_{[n]}^+ (u^{(\Lambda)} ) \ts \Ec_\Lambda ,\label{lemma61a1} 
\intertext{\item The following relations   hold  for the $R$-matrices \eqref{erofzed} and \eqref{erofu}:}
&R_{1n}^{12}(z e^{-u^{(\Lambda)}+dh})^{\pm 1} \ts \Ec_\Lambda =\Ec_\Lambda\ts  R_{1n}^{12}(z e^{-u^{(\Lambda)}+dh})^{\pm 1}  \ts \Ec_\Lambda  \quad\text{for all }d\in\CC,\label{lemma61c}\\
&R_{1n}^{12}(  e^{v-u^{(\Lambda)}+dh})^{\pm 1} \ts \Ec_\Lambda =\Ec_\Lambda\ts  R_{1n}^{12}(  e^{v-u^{(\Lambda)}+dh})^{\pm 1}  \ts \Ec_\Lambda  \quad\text{for all }d\in\CC,\label{lemma61c0}
\intertext{where $\Ec_\Lambda$ is applied on the tensor factors $2,\ldots ,n+1$ of   $\ndo\CC^N \ot(\ndo\CC^N)^{\ot n}$.}
\intertext{\item The following relation    holds for the $R$-matrix \eqref{erofzed}:}
&R_{n1}^{12}(z e^{ u^{(\Lambda)}+dh})^{\pm 1} \ts \Ec_\Lambda =\Ec_\Lambda\ts  R_{n1}^{12}(z e^{ u^{(\Lambda)}+dh})^{\pm 1}  \ts \Ec_\Lambda \quad\text{for all }d\in\CC,
\label{lemma61d}
\end{align}
 where $\Ec_\Lambda$ is applied on the tensor factors $1,\ldots ,n$ of   $ (\ndo\CC^N)^{\ot n}\ot \ndo\CC^N$. 
\end{enumerate}
\end{lem}

Suppose   $\Gamma$ is a  standard tableau of shape $\gamma\vdash m,$ and let $\Ec_{\Lambda,\Gamma}=\Ec_{\Lambda}\ot \Ec_{ \Gamma}.$ By combining \eqref{lemma61c} and  \eqref{lemma61d}, one obtains the following lemma.

\begin{lem}\label{Lemma_62}
The following relation    holds for the $R$-matrix \eqref{erofzed}:
$$
R_{nm}^{12}(z e^{ u^{(\Lambda)}-v^{(\Gamma)}+dh})^{\pm 1} \ts \Ec_{\Lambda,\Gamma} =\Ec_{\Lambda,\Gamma}\ts  R_{nm}^{12}(z e^{u^{(\Lambda)}-v^{(\Gamma)}+dh})^{\pm 1}  \ts \Ec_{\Lambda,\Gamma} \quad\text{for all }d\in\CC.
$$
\end{lem}

 Motivated by the construction of commutative families in the  $q$-Yangian of type $A$   \cite[Cor. 3.4]{JLM}, we introduce the formal power series
\beq\label{telambda}
T_\Lambda(u)=\tr_{1,\ldots ,n} \, T_{[n]}^+ (u^{(\Lambda)} )\vac \ts D^{\ot n}\ts \Ec_\Lambda \in \mathcal{V}^{crit}(R)[[u]].
\eeq
In addition, we write $T_{\Lambda,\Gamma}(u,v)=T_{\Lambda}(u)\ot T_{\Gamma }(v) \in \mathcal{V}^{crit}(R)^{\ot 2}[[u,v]] .$
The next proposition  is  a trigonometric counterpart of \cite[Lemma 4.9]{BK0}.

\begin{pro}\label{prop63}
The maps $ \sigma$ and  $\rho ,$ as given by Proposition \ref{lemma_31}, satisfy
$$
\sigma(z)\left( T_{\Lambda,\Gamma}(u,v)\right)
=\rho(z)\left( T_{\Lambda,\Gamma}(u,v)\right)
=T_{\Lambda,\Gamma}(u,v).
$$
\end{pro}

\begin{prf}
Let us prove 
$$\sigma(z)\left( T_{\Lambda,\Gamma}(u,v)\right)
=T_{\Lambda,\Gamma}(u,v).$$
 By \eqref{sigma_formula}, we have
\begin{align}
\sigma(z)\left( T_{\Lambda,\Gamma}(u,v)\right)
=\tr_{1,\ldots ,n+m}\ts R_{nm}(z,u,v)\ts T_{nm}(u,v)\ts R_{nm}(z,u,v)^{-1}\ts D^{\ot (n+m)}\ts\Ec_{\Lambda,\Gamma},\label{rhs_sigma_1}
\end{align}
where
\begin{align}
&R_{nm}(z,u,v)^{\pm 1} = R_{nm}^{12}(ze^{u^{(\Lambda)} -v^{(\Gamma)}})^{\pm 1},\label{rhs_sigma_2} \\
 &T_{nm}(u,v)=  T_{1\ts n+m+1}^+ (u-c_1(\Lambda)h )\ldots T_{n\ts n+m+1}^+ (u-c_n(\Lambda)h )\non\\
 &\qquad\qquad\quad\times T_{n+1\ts n+m+2}^+ (v-c_1(\Gamma)h )\ldots T_{n+m\ts n+m+2}^+ (v-c_m(\Gamma)h ).\non
\end{align}
We now rewrite the right-hand side of \eqref{rhs_sigma_1} as follows. First, we use 
$R_{12}(z)  D_1 D_2 =D_1 D_2 R_{12}(z)   $ to move the term $D^{\ot (n+m)}$ to the left, thus getting
$$
\tr_{1,\ldots ,n+m}\ts R_{nm}(z,u,v)\ts T_{nm}(u,v)\ts D^{\ot (n+m)}\ts R_{nm}(z,u,v)^{-1}\ts\Ec_{\Lambda,\Gamma}.
$$
By Lemma \ref{Lemma_62}, this equals
$$
\tr_{1,\ldots ,n+m}\ts R_{nm}(z,u,v)\ts T_{nm}(u,v)\ts D^{\ot (n+m)}\ts\Ec_{\Lambda,\Gamma}\ts R_{nm}(z,u,v)^{-1}\ts\Ec_{\Lambda,\Gamma}.
$$
Due to the cyclic property of the trace, we can move $\Ec_{\Lambda,\Gamma}$ to the left, so that we obtain
$$
\tr_{1,\ldots ,n+m}\ts \Ec_{\Lambda,\Gamma}\ts R_{nm}(z,u,v)\ts T_{nm}(u,v)\ts D^{\ot (n+m)}\ts\Ec_{\Lambda,\Gamma}\ts R_{nm}(z,u,v)^{-1}.
$$
This is equal to
$$
\tr_{1,\ldots ,n+m}\ts \Ec_{\Lambda,\Gamma}\ts R_{nm}(z,u,v)\ts T_{nm}(u,v)\ts \Ec_{\Lambda,\Gamma}\ts  D^{\ot (n+m)}\ts R_{nm}(z,u,v)^{-1},
$$
since $ \Ec_{\Lambda,\Gamma}$ and $D^{\ot (n+m)}$ commute. Next, we apply \eqref{lemma61a1}:
$$
\tr_{1,\ldots ,n+m}\ts \Ec_{\Lambda,\Gamma}\ts R_{nm}(z,u,v)\ts \Ec_{\Lambda,\Gamma}\ts T_{nm}(u,v)\ts \Ec_{\Lambda,\Gamma}\ts  D^{\ot (n+m)}\ts R_{nm}(z,u,v)^{-1}.
$$
By employing Lemma \ref{Lemma_62} once again, we get
$$
\tr_{1,\ldots ,n+m}\ts   R_{nm}(z,u,v)\ts \Ec_{\Lambda,\Gamma}\ts T_{nm}(u,v)\ts \Ec_{\Lambda,\Gamma}\ts  D^{\ot (n+m)}\ts R_{nm}(z,u,v)^{-1}.
$$
Using the cyclic property of the trace, we move the term $R_{nm}(z,u,v)$ to the right and then we cancel the corresponding $R$-matrices, so that we get
\begin{align*}
=&\,\,\tr_{1,\ldots ,n+m}\ts \Ec_{\Lambda,\Gamma}\ts T_{nm}(u,v)\ts \Ec_{\Lambda,\Gamma}\ts  D^{\ot (n+m)}\ts R_{nm}(z,u,v)^{-1}\ts   R_{nm}(z,u,v) \\
=&\,\,\tr_{1,\ldots ,n+m}\ts \Ec_{\Lambda,\Gamma}\ts T_{nm}(u,v)\ts \Ec_{\Lambda,\Gamma}\ts  D^{\ot (n+m)} .
\end{align*}
Clearly, this equals
$$
\tr_{1,\ldots ,n+m}\ts  T_{nm}(u,v)\ts \Ec_{\Lambda,\Gamma}\ts  D^{\ot (n+m)}\ts \Ec_{\Lambda,\Gamma} .
$$
Finally,  by using the fact that $ \Ec_{\Lambda,\Gamma}$ and $D^{\ot (n+m)}$ commute and then the identity 
 $\left( \Ec_{\Lambda,\Gamma}\right)^2 = \Ec_{\Lambda,\Gamma},$ the above expression turns to   
$$
\tr_{1,\ldots ,n+m}\ts  T_{nm}(u,v) \ts  D^{\ot (n+m)} \ts \Ec_{\Lambda,\Gamma}= T_{\Lambda,\Gamma}(u,v),
$$
as required.

Let us verify the second assertion of the proposition,
\beq\label{secondassertion}
\rho(z)\left( T_{\Lambda,\Gamma}(u,v)\right)
=T_{\Lambda,\Gamma}(u,v).
\eeq
 By \eqref{rho_formula}, the expression $\rho(z)\left( T_{\Lambda,\Gamma}(u,v)\right)$ is equal to
 \begin{align}
\tr_{1,\ldots ,n+m}\ts T_{\Lambda}^{13}(u)\ts R_{nm}(ze^{-Nh},u,v)^{-1}\ts T_{\Gamma}^{24}(v)\ts R_{nm}(z,u,v)\ts D^{\ot (n+m)}\ts\Ec_{\Lambda,\Gamma},\label{rhs_rho_1}
\end{align}
where we use the notation \eqref{rhs_sigma_2} along with
\begin{align*}
 T^{13}_{\Lambda}(u )=&\,\ts T_{1\ts n+m+1}^+ (u-c_1(\Lambda)h )\ldots T_{n\ts n+m+1}^+ (u-c_n(\Lambda)h ),\\
  T^{24}_{\Gamma}(v )=&\, T_{n+1\ts n+m+2}^+ (v-c_1(\Gamma)h )\ldots T_{n+m\ts n+m+2}^+ (v-c_m(\Gamma)h ).
\end{align*}
To prove \eqref{secondassertion}, it suffices to show that  \eqref{rhs_rho_1} equals
\beq\label{rhs_rho_2}
\tr_{1,\ldots ,n+m}\ts T_{\Lambda}^{13}(u)\ts\Ec_{\Lambda }\ts  R_{nm}(ze^{-Nh},u,v)^{-1}\ts\Ec_{ \Gamma}\ts T_{\Gamma}^{24}(v)\ts R_{nm}(z,u,v)\ts D^{\ot (n+m)},
\eeq
where $\Ec_{\Lambda }$ (resp. $\Ec_{ \Gamma}$) is applied on the tensor factors $1,\ldots ,n$ (resp. $n+1,\ldots ,n+m$). Indeed, \eqref{rhs_rho_2} coincides with
\beq\label{rhs_rho_3}
\tr_{1,\ldots ,n+m}\ts T_{\Lambda}^{13}(u)\ts\Ec_{\Lambda , \Gamma}\ts T_{\Gamma}^{24}(v)\ts \left(R_{nm}(ze^{-Nh},u,v)^{-1}\cdotlr\left( R_{nm}(z,u,v)\ts D^{\ot (n+m)}\right)\right),
\eeq
where the symbol ``$\cdotlr$'' denotes the standard multiplication in the algebra
$(\ndo\CC^N)^{\ot n}\ot ((\ndo\CC^N)^{op})^{\ot m} $ and $(\ndo\CC^N)^{op}$ stands for the opposite algebra of $\ndo\CC^N .$ By the crossing symmetry property \eqref{csym}, we have
$$
R_{nm}(ze^{-Nh},u,v)^{-1}\cdotlr\left( R_{nm}(z,u,v)\ts D^{\ot (n+m)}\right)=D^{\ot (n+m)},
$$
so that
\eqref{rhs_rho_3} is equal to
$$
\tr_{1,\ldots ,n+m}\ts T_{\Lambda}^{13}(u)\ts\Ec_{\Lambda , \Gamma}\ts T_{\Gamma}^{24}(v) \ts D^{\ot (n+m)}
= \tr_{1,\ldots ,n+m}\ts  T_{nm}(u,v) \ts D^{\ot (n+m)}\ts\Ec_{\Lambda , \Gamma}=T_{\Lambda,\Gamma}(u,v),
$$
as required.
Hence, to finish the proof, it remains to check that    the expressions  \eqref{rhs_rho_1} and  
\eqref{rhs_rho_2} coincide. This follows by the arguments which go in parallel with the  proof of the first assertion of the theorem and rely on Lemmas \ref{lemma_61} and \ref{Lemma_62}.
\end{prf}
 
 Proposition \ref{prop63}, together with \eqref{shtsm} and \eqref{s_hat_identity}, implies
that the map \eqref{Shatmaptrig} satisfies
$$
\wht{\Sc}(z)\left( T_{\Gamma,\Lambda}(u,v)\right)
=T_{\Gamma,\Lambda}(u,v),
$$
so that, for the braiding \eqref{braiding} of $\mathcal{V}^{crit}(R)$, we also have
	$$
 \Sc (z)\left( T_{\Gamma,\Lambda}(u,v)\right)
=T_{\Gamma,\Lambda}(u,v).
$$
 Hence, the vertex operators of coefficients of all series $T_{ \Lambda}(u )$ are mutually local in the $h$-adic sense. In the next subsection, we shall further strengthen this result.

 %%%%%%%%%%%%%%%%%%%%%%%%%%%%%%%%%%%
\subsection{Central elements in  the  \texorpdfstring{$h$}{h}-adic quantum vertex algebra \texorpdfstring{$\mathcal{V}^{crit}(R^{trig})$ }{Vcrit(Rtrig)}}\label{sec0603}
%%%%%%%%%%%%%%%%%%%%%%%%%%%%%%%%%%%

In this subsection, we consider the {\em center} of the $h$-adic quantum vertex algebra 
$\mathcal{V}^{crit}(R ),$
$$
\mathfrak{z}(\mathcal{V}^{crit}(R ))=
\left\{v\in \mathcal{V}^{crit}(R )\,:\, Y(u,z)v\in \mathcal{V}^{crit}(R )[[z]] \text{ for all }u\in \mathcal{V}^{crit}(R )\right\}.
$$
For more information on the notion of center   see \cite[Sect. 4]{DGK} and \cite[Sect. 3.2]{JKMY}. 
The next proposition generalizes \cite[Prop. 3.5]{KM}, where the families of central elements were constructed in the case when  $\Lambda  $ consists of a single column, i.e., when $\Ec_\Lambda$ is the action of the normalized anti-symmetrizer.

\begin{pro}\label{prop_64}
All coefficients of the series $T_{\Lambda}(u)$ belong to the center of $\mathcal{V}^{crit}(R ) .$
\end{pro}

\begin{prf}
To show that an element $w\in\mathcal{V}^{crit}(R )$ belongs to the center $\mathfrak{z}(\mathcal{V}^{crit}(R )),$ it suffices to check that it satisfies
$T^-(u)w=w.$ Indeed, this immediately follows from the expression \eqref{qva1} for the vertex operator map of $\mathcal{V}^{crit}(R )$. Hence, to prove the proposition, it is sufficient to verify the  identity 
$$
T^-(u)\ts T_\Lambda(v) = T_\Lambda(v) .
$$

By using the definitions    \eqref{teminus} and \eqref{telambda} of $T^-(u)$ and  $ T_\Lambda(v)$, we find
\begin{align}
T^-(u)\ts T_\Lambda(v) =&\,\,
 \tr_{1,\ldots ,n}\, T_0^-(u)\ts T_{[n]}^+ (v^{(\Lambda)} )\vac \ts D^{\ot n}\ts \Ec_\Lambda \non\\
=&\,\,  \tr_{1,\ldots ,n}\, R_{1n}^{01}(e^{u-v(\Lambda)-Nh/2})^{-1}\ts T_{[n]}^+ (v^{(\Lambda)} )\vac \ts R_{1n}^{01}(e^{u-v(\Lambda)+Nh/2}) \ts D^{\ot n}\ts \Ec_\Lambda .\label{prop_64_eq1}
\end{align}
Next, by arguing as in the proof of the first assertion of Proposition \ref{prop63}, one can show that \eqref{prop_64_eq1} is equal to
\beq\label{prop_64_eq2}
\tr_{1,\ldots ,n}\, R_{1n}^{01}(e^{u-v(\Lambda)-Nh/2})^{-1}\ts T_{[n]}^+ (v^{(\Lambda)} )\vac  \ts D^{\ot n}\ts \Ec_\Lambda \ts R_{1n}^{01}(e^{u-v(\Lambda)+Nh/2}).
\eeq
Naturally, instead of Lemma \ref{Lemma_62}, the argument now employs the identity \eqref{lemma61d}, as we are   using the (additive) $R$-matrix \eqref{erofu}.
Finally, to conclude that \eqref{prop_64_eq2} coincides with
$$
\tr_{1,\ldots ,n}\,   T_{[n]}^+ (v^{(\Lambda)} )\vac  \ts D^{\ot n}\ts \Ec_\Lambda= T_\Lambda(v)   ,
$$
i.e., to cancel the $R$-matrix terms in \eqref{prop_64_eq2},
one uses the crossing symmetry property \eqref{csym_reu} in the same way as in the proof of   the second assertion of Proposition \ref{prop63}.
\end{prf}

%%%%%%%%%%%%%%%%%%%%%%%%%%%%%%%%%%%
\subsection{Commutative families in the \texorpdfstring{$h$}{h}-Yangian  \texorpdfstring{$\Y$}{Yh(glN)}}\label{sec0604}
%%%%%%%%%%%%%%%%%%%%%%%%%%%%%%%%%%%

In this subsection, we demonstrate an application of Theorem \ref{thm_53-trigg} and Proposition \ref{prop63} to constructing commutative families in the $h$-Yangian. Roughly speaking, the  commutative families emerge as images of families of central elements of $\mathcal{V}^{crit}(R )  ,$ established by Proposition \ref{prop_64}, under certain $(\sigma   ,\rho    )$-deformed $\phi$-coordinated module maps.

Let $I^{(p)}$ with $p\geqslant 1$ be the $h$-adically completed ideal in the algebra $\Y$ generated by all
$\ell_{i j}^{(r)}$ with $r>p.$   The quotient $W^{(p)}\coloneqq \Y / I^{(p)}$ is naturally equipped with the structure of restricted $\Y$-module.  Therefore, by Theorem \ref{thm_53-trigg}, there exists a  structure $(W^{(p)},Y_{W^{(p)}})$ of $(\sigma   ,\rho    )$-deformed $\phi$-coordinated $\mathcal{V}^{crit}(R^{trig})$-module over $W^{(p)},$ which is uniquely determined by
 \eqref{thm53eq}. In particular, we have
$$
Y_{W^{(p)}}(T_\Lambda(0),z)
=\tr_{1,\ldots ,m}\, L_{[m]}^- (z^{(\Lambda)} )  \ts D^{\ot m}\ts \Ec_\Lambda \quad\text{for all }p\geqslant 1.
$$
Set
\beq\label{eloflambdda}
L_\Lambda(z)_{W^{(p)}} = Y_{W^{(p)}}(T_\Lambda(0),z).
\eeq
By the $\sigma$-locality \eqref{sigma-loc-trigg} and Proposition \ref{prop63}, we have
\beq\label{eloflambddaa}
L_\Lambda(z_1)_{W^{(p)}}\ts L_\Gamma(z_2)_{W^{(p)}}
=L_\Gamma(z_2)_{W^{(p)}}\ts L_\Lambda(z_1)_{W^{(p)}}\quad\text{for all }p\geqslant 1.
\eeq
Indeed, the powers of $(z_1-z_2),$ which come from the  $\sigma$-locality, cancel, as the expression in
\eqref{eloflambdda}
is a polynomial in $z^{-1} .$ Hence, since the intersection $\cap_{p\geqslant 1} I^{(p)}$ is trivial, we conclude from \eqref{eloflambddaa} that the identity
$$
L_\Lambda(z_1) \ts L_\Gamma(z_2) 
=L_\Gamma(z_2) \ts L_\Lambda(z_1) 
$$
holds for
$$
L_\Lambda(z)= \tr_{1,\ldots ,m}\, L_{[m]}^- (z^{(\Lambda)} )  \ts D^{\ot m}\ts \Ec_\Lambda \in \Y[[z^{-1}]]
.
$$
Thus, we proved the next proposition, which is a Yangian counterpart of \cite[Cor. 3.4]{JLM}.
\begin{pro}
The coefficients of all series $L_\Lambda(z)$ pairwise commute.
\end{pro}

%%%%%%%%%%%%%%%%%%%%%%%%%%%%%%%%%%%

\section*{Acknowledgement} 
L. B. is member of Gruppo Nazionale per le Strutture Algebriche, Geometriche e le loro Applicazioni  (GNSAGA) of the Istituto Nazionale di Alta Matematica (INdAM).
The research reported in this paper was finalized during the third author’s visit to the School of Mathematics and Statistics at the Central China Normal
University in Wuhan. He is grateful to the School for
warm hospitality. 
S. K. is partially supported by the Croatian Science
Foundation under the project IP-2025-02-4720  and by the project ``Implementation of cutting-edge research and its application as part of the Scientific Center of Excellence for Quantum and Complex Systems, and Representations of Lie Algebras'', Grant No. PK.1.1.10.0004, co-financed by the European Union through the European Regional Development Fund - Competitiveness and Cohesion Programme 2021--2027.

 \linespread{1.0}

\end{document}